\documentclass[12pt]{article}
\usepackage{amsfonts}

\usepackage{latexsym}
\usepackage{amssymb}
\usepackage{epsfig}
\usepackage{amsmath}
\usepackage{amsthm}
\usepackage[mathscr]{eucal}
\usepackage{subfigure}
\usepackage{graphicx}

\newtheorem{Theorem}{Theorem}
\newtheorem{Proposition}{Proposition}
\newtheorem{Lemma}{Lemma}[section]
\newtheorem{Corollary}{Corollary}

\renewcommand{\qed}{\hfill{\ \ \rule{2mm}{2mm}} \vspace{0.2in}}

\newcommand{\ind}{1\hspace{-2.3mm}{1}}

\begin{document}

\title{Convergence in First Passage Percolation with nonidentical passage times}
\author{ \textbf{Ghurumuruhan Ganesan}
\thanks{E-Mail: \texttt{gganesan82@gmail.com}} \\
%EndAName
\ \\
New York University, Abu Dhabi}
\date{}
\maketitle

\begin{abstract}
In this paper we consider first passage percolation on the square lattice \(\mathbb{Z}^d\) with
edge passage times that are independent and have uniformly bounded second moment, but not necessarily identically distributed. For integer \(n \geq 1,\) let \(T_n\) be the minimum passage time between the origin and the point \((n,0,\ldots,0).\) We prove that \(\frac{1}{n}(T_n-\mathbb{E}T_n)\) converges to zero almost surely and in \(L^2\) as \(n~\rightarrow~\infty.\) The convergence is nontrivial in the sense that \(\frac{T_n}{n}\) is asymptotically bounded away from zero and infinity almost surely. We first define a truncated version \(\hat{T}^{(n)}_n\) that is asymptotically equivalent to~\(T_n.\) We then use a finite box modification of the martingale method of Kesten~(1993) to estimate the variance of \(\hat{T}^{(n)}_n.\) Finally, we use a subsequence argument to obtain almost sure convergence for \(\frac{1}{n}(\hat{T}^{(n)}_n - \mathbb{E}\hat{T}^{(n)}_n).\) The corresponding result for \(T_n\) is then obtained using asymptotic equivalence of \(T_n\) and \(\hat{T}^{(n)}_n.\) For identically distributed passage times, our method alternately obtains almost sure convergence of~\(\frac{T_n}{n}\) to a positive constant~\(\mu_F,\) without invoking the subadditive ergodic theorem.

\vspace{0.1in} \noindent \textbf{Key words:} First passage percolation,
nonidentical passage times, almost sure convergence, subsequence argument.

\vspace{0.1in} \noindent \textbf{AMS 2000 Subject Classification:} Primary:
60J10, 60K35; Secondary: 60C05, 62E10, 90B15, 91D30.
\end{abstract}

\bigskip

\newtheorem{thm}{Lemma}[section]
\newtheorem*{cor}{Corollary}

\setcounter{equation}{0}
\renewcommand\theequation{\thesection.\arabic{equation}}
\section{Introduction} \label{intro}
Consider the square lattice \(\mathbb{Z}^d,\) where two vertices \(w_1 = (w_{1,1},\ldots,w_{1,d})\) and \(w_2 = (w_{2,1},\ldots,w_{2,d})\) are adjacent if \(\sum_{i=1}^{d}|w_{1,i} - w_{2,i}| = 1\) and adjacent vertices are joined together by an edge. Let \(\{q_i\}_{i \geq 1}\) denote the set of edges. Each edge \(q_i\) is equipped with a random passage time \(t(q_i)\) and for integer \(n \geq 1,\) we are interested in the shortest passage time \(T_n\) from the origin to the point \((n,0,\ldots,0).\) We give formal definitions below.

When the passage times are independent and identically distributed (i.i.d.) the subadditive ergodic theorem (Kingman (1973)) is used for studying almost sure convergence and convergence in mean of \(\frac{T_n}{n}\) (see Smythe and Wierman~(2008), Cox and Durrett (1981), Kesten (1986) and references therein). In many cases of interest, it may happen that the passage times are not i.i.d. As a simple example we think of the nodes as mobile stations and passage time as the time taken to send a packet from one node to another. Links between adjacent nodes may have different passage time depending on the geographical conditions~etc. In such cases, it is of interest to study convergence properties of the first passage time \(T_n\) with appropriate centering and scaling.

%In our main result (Theorem~\ref{thm1}), we show that even if the passage times are not i.i.d., the centred sequence \(\frac{1}{k}(T_k-\mathbb{E}T_k)\) converges to zero  in \(L^2\) and almost surely, provided the passage times satisfy some mild moment conditions. As a corollary (Corollary~\ref{cor1}), for the i.i.d. case, we obtain the almost sure convergence of \(\frac{T_k}{k}\) to a constant \(\mu,\) without invoking the subadditive ergodic theorem.

%For the edge \(q_i\) of \(\mathbb{Z}^d,\) the random passage time \(t(q_i)\) is nonnegative and let \(\mathbb{P}_i\) probability space \((\Omega_i, {\cal B}_i, \mathbb{P}_i).\)

\subsection*{Model}
We briefly describe the probability space first. For integer \(i \geq 1,\) let \(\Omega_i = \mathbb{R}\) and \({\cal B}_i = {\cal B}(\Omega_i)\) denote real line and the Borel sigma field, respectively. Let \(\Omega = \otimes_{i=1}^{\infty} \Omega_i\) and \({\cal F} = \otimes_{i=1}^{\infty} {\cal B}(\Omega_i)\) denote the product space and product sigma field, respectively. Here \({\cal F}\) is the product sigma algebra generated by the cylinder sets of the form \(\otimes_{i=1}^{\infty} A_i\) where each~\(A_i\) is a Borel set in \(\mathbb{R}\) and \(A_i = \mathbb{R}\) for all but a finite set of values of~\(i.\) We recall that \(t(q_i)\) is the random passage time of the edge \(q_i\) of \(\mathbb{Z}^d.\) We define the random sequence \((t(q_1),t(q_2),\ldots)\) on the probability space \((\Omega, {\cal F}, \mathbb{P}).\) If \(\omega \in \Omega\) is a realization of the passage times \((t(q_1),t(q_2),\ldots),\) we say that \(t(q_i) = t(q_i,\omega)\) is the passage time of the edge~\(q_i\) for the realization \(\omega.\)

In what follows, we consider passage times of paths and we therefore give a formal definition below. A \emph{path} \(\pi\) in \(\mathbb{Z}^d\) is a sequence of distinct edges \((e_1,...,e_t)\) in \(\mathbb{Z}^d\) with the following three properties: The edge \(e_1\) shares an endvertex only with edge \(e_2\) and no other edge in \(\pi.\) The edge \(e_t\) shares an endvertex with only the edge \(e_{t-1}\) and no other edge in \(\pi.\) For \(2 \leq i \leq t-1,\) the edge~\(e_i\) shares an endvertex with only the edges \(e_{i-1}\) and \(e_{i+1}\) and no other edge in~\(\pi.\) All paths we consider in this paper are self avoiding paths with finite number of edges. Further let \(a\) be the endvertex of \(e_1\) not common with \(e_2\) and let \(b\) be the endvertex of \(e_t\) not common with~\(e_{t-1}.\) We say that \(a\) and \(b\) are the \emph{endvertices} of the path \(\pi.\)

For \(\omega \in \Omega\) and a path \(\pi = (e_1,\ldots,e_r), e_i \subset \{q_j\}\) containing \(r\) edges, we define the passage time of~\(\pi\) as
\begin{equation}\label{t_pi_def}
T(\pi) = T(\pi, \omega) = \sum_{i=1}^{r} t(e_i, \omega).
\end{equation}
Letting \({\cal P}_n\) denote the set of all finite paths with endvertices origin and \((n,0,...,0),\) we define
\begin{equation}\label{t_0n_def}
T_n = T_n(\omega) = \inf_{\pi \in {\cal P}_n} T(\pi,\omega)
\end{equation}
to be the \emph{minimum passage time} between the origin and \((n,0,\ldots,0).\) For convenience, we suppress the dependence on~\(\omega\) unless specifically mentioned.

% not a cylinder set....
To see the measurability of \(T_n,\) we fix a finite path \(\pi\) and let \(T(\pi)\) be the corresponding passage time as defined in (\ref{t_pi_def}) above. For finite paths \(\pi\) the passage time \(T(\pi) = \sum_{e \in \pi} t(e)\) is simply the sum of the passage times of the individual edges. Since each~\(t(e)\) is \({\cal F}-\)measurable, we have that \(T(\pi)\) is \({\cal F}-\)measurable. Also, we have that \({\cal P}_n = \cup_{m \geq 1} {\cal P}_{n,m}\) where~\({\cal P}_{n,m}\) denotes the set of all paths contained in the box~\(B_m = [-m,m]^d\) and having endvertices as the origin and the point \((n,0,\ldots,0).\) Therefore \({\cal P}_n\) is countable and we have from~(\ref{t_0n_def}) that \(T_n\) is also \({\cal F}-\)measurable.

%We have that \(T_n\) defined in (\ref{t_0n_def}) is \(\otimes_{i=1}^{\infty} {\cal F}_i\) measurable. For completeness, we give a small proof. We note that the passage time \(T_n\) can also be written as the limit \(T_n = \lim_{M \rightarrow \infty} T_M(n)\) where \(T_m(n) := \min_{\pi \in {\cal P}_{n,m}} T(\pi)\) and \({\cal P}_{n,m}\) denotes the set of all paths starting from origin and ending at \((n,0,\ldots,0)\) and contained in \(B_n.\) Since \(T_M(n)\) is \(\otimes_{i=1}^{\infty} {\cal F}_i\) measurable for each integer \(M \geq 1,\) we have that \(T_n\) is also \(\otimes_{i=1}^{\infty} {\cal F}_i\) measurable. % and we define the first passage percolation process on the probability spac.

Our aim is to study convergence (almost surely and in mean) of the sequence \(T_n\) suitably centred and scaled, when the passage times \(\{t(q_i)\}\) are independent (but not necessarily identically distributed) random variables. We define the following mild conditions on the passage times.\\
\((i)\) We have that \(\sup_{i} \mathbb{P}(t(q_i) < \epsilon)  \longrightarrow 0\) as \(\epsilon \downarrow 0.\)\\
\((ii)\) We have that \(\sup_{i} \mathbb{E}t^2(q_i) < \infty.\)\\
\((ii)(a)\) We have that \(\{t^2(q_i)\}_{i \geq 1}\) is a uniformly integrable sequence in the sense that \[\sup_{i \geq 1} \mathbb{E}t^2(q_i)\ind(t(q_i) \geq M) \longrightarrow 0\] as \(M \rightarrow \infty.\)\\
\((ii)(b)\) We have that \(\sup_{i} \mathbb{E}t^p(q_i) < \infty\) for some \(p > 2.\)\\\\

The conditions \((ii)(a)\) and \((ii)(b)\) are needed only for the \(L^2-\)convergence results and are stronger than condition \((ii);\) i.e., condition \((ii)(b)\) implies condition \((ii)(a)\) implies condition \((ii).\) Unless otherwise mentioned, all results in this paper are derived assuming only conditions \((i)-(ii).\)
%In Section~\ref{esti} below, we use conditions \((i)-(ii)\) to show that the infimum in (\ref{t_0n_def}) is almost surely a minimum.

%\begin{array}{clrr}%
%a+b+c & uv & x-y & 27 \\       x+y  & w  & +z  & 363
%\end{array}

The following is the main result of our paper.
\begin{Theorem} \label{thm1} If conditions \((i)\) and \((ii)\) hold, we have that
\begin{equation} \label{conv_tn_main}
\frac{1}{n}\left(T_n - \mathbb{E}T_n\right) \longrightarrow 0\;\;\text{a.s. and in }L^1
\end{equation}
as \(n \rightarrow \infty.\)  If conditions \((i)\) and \((ii)(a)\) hold, then (\ref{conv_tn_main}) holds a.s. and in \(L^2.\) If conditions \((i)\) and \((ii)(b)\) hold, then (\ref{conv_tn_main}) holds a.s. and in \(L^2\) and also
\begin{equation}\label{var_lin_11}
var(T_n-\mathbb{E}T_n) \leq C n
\end{equation}
for some constant \(C > 0\) and for all \(n \geq 1.\) In all the cases, the convergence is nontrivial in the sense that there are constants \(\eta_1,\eta_2 \in (0,\infty) \) such that
\begin{equation}\label{tmin}
\eta_1 \leq \liminf_n \frac{T_n}{n}  \leq \limsup_n \frac{T_n}{n} \leq \eta_2 \text{ a.s.}
\end{equation}
\end{Theorem}

For the particular case of i.i.d. passage times, the uniform integrability condition \((ii)(a)\) is implied by the moment condition \((ii)\) and we have the following as a Corollary of Theorem~\ref{thm1}.
\begin{Corollary} \label{cor1} If the passage times are i.i.d. and conditions \((i)\) and \((ii)\) hold, we have that
\begin{equation}
\frac{T_n}{n}  \longrightarrow \mu_F\;\;\text{a.s. and in }L^2
\end{equation}
as \(n \rightarrow \infty,\) for some constant \(\mu_F > 0.\)
\end{Corollary}
The constant \(\mu_F\) is also called the time constant. We remark that an important contribution in our paper is the use of truncation and a subsequence argument (described below) to obtain almost sure convergence for \(\frac{1}{n}(T_n - \mathbb{E}T_n).\) For the particular case of i.i.d. passage times, we alternately obtain almost sure convergence of~\(\frac{T_n}{n}\) to the constant~\(\mu_F,\) without invoking the subadditive ergodic theorem. For more material on first passage percolation using subadditive ergodic theorem, we refer to Alexander (1993), Cox and Durrett (1981), Kesten~(1993) and Smythe and Wierman (2008) and references therein.

We give a brief outline of the idea of the proof presented in Section~\ref{pf1}.
For a fixed \(n \geq 1,\) we first truncate the passage time of every edge by~\(n^{\alpha}\) for some constant \(\alpha > 0.\) For appropriately chosen \(\alpha,\) we show that  the corresponding truncated minimum passage time~\(\hat{T}^{(n)}_n\) is asymptotically equivalent to~\(T_n.\) Using a finite box modification of the martingale difference method of Kesten~(1993), we then show that the variance of \(\hat{T}^{(n)}_n\) grows at most linearly with \(n.\) We use the variance estimate and a subsequence argument to show that \(\frac{\hat{T}^{(n)}_n-\mathbb{E}\hat{T}^{(n)}_n}{n} \longrightarrow 0\) a.s. and in \(L^2\) as \(n \rightarrow \infty.\) From the asymptotic equivalence of \(T_n\) and \(\hat{T}^{(n)}_n,\) we then obtain the analogous result for~\(T_n.\)

%//CHECK ALL + ETC.. .AGN AND AGN...

%For \(n \geq 1,\) we first truncate the passage time of every edge by \(n^{\alpha}\) and for \(\alpha > 0\) appropriately chosen, we show that the truncated minimum passage time \(\hat{T}^{(n)}_n\) is asymptotically equivalent to~\(T_n.\) We then use a finite box modification of the martingale difference method of Kesten~(1993) to obtain a variance estimate for~\(\hat{T}^{(n)}_n.\) Finally, we use the variance estimate and a subsequence argument to show that \(\frac{1}{n}(\hat{T}^{(n)}_n-\mathbb{E}\hat{T}^{(n)}_n)\) converges to zero almost surely and in \(L^2\) as \(n \rightarrow \infty.\) The analogous result for~\(T_n\) then follows from the asymptotic equivalence of \(T_n\) and~\(\hat{T}^{(n)}_n.\)

The paper is organized as follows: In Section~\ref{esti}, we define the construction of the auxiliary truncated passage time \(\hat{T}^{(n)}_n\) that is useful to study the convergence of \(T_n.\) We also obtain auxiliary results needed for future use. In Section~\ref{pf1}, we prove that the truncated and untruncated passage times are asymptotically equivalent and that it suffices to study the a.s. convergence of \(\frac{1}{n}\left(\hat{T}^{(n)}_n - \mathbb{E}\hat{T}^{(n)}_n\right).\) In Section~\ref{var_est_sec} we obtain variance estimates for both \(\hat{T}^{(n)}_n\) and~\(T_n,\) to be used in the proof of Theorem~\ref{thm1}. Finally in Section~\ref{pf_main_th}, we prove Theorem~\ref{thm1} and  Corollary~\ref{cor1}. %We conclude in Section~\ref{future} with a brief discussion on possible future directions.

\setcounter{equation}{0}
\renewcommand\theequation{\thesection.\arabic{equation}}
\section{Truncated passage time and geodesic estimates}~\label{esti}
Throughout the paper, we use an approximate truncated passage time~\(\hat{T}^{(n)}_n\) to obtain estimates on \(T_n.\) In this section, we define~\(\hat{T}^{(n)}_n\) formally and show that almost surely minimum passage times are attained by finite paths (geodesics). We also obtain estimates on the lengths of geodesics that are useful in the proof of Theorem~\ref{thm1}.

%obtain estimates on \(var(\hat{T}^{(n)}_n)\) and \(var(T_n)\) used in the proof of Theorem~\ref{thm1}

For integer \(k \geq 1,\) we define auxiliary random variables \(\{\hat{T}^{(k)}_n\}_{n \geq 1}\) as follows. As before let \(\omega \in \Omega\) be any fixed realization of the passage times \((t(q_1),t(q_2),\ldots)\) of the edges \(\{q_i\}\) of \(\mathbb{Z}^d.\) For \(i \geq 1,\) let
\begin{equation}\label{trunc_pas}
t^{(k)}(q_i) = t^{(k)}(q_i,\omega) = \min(t(q_i,\omega), k^{\alpha}),
\end{equation}
where \(\alpha > 0\) is a constant to be determined later and as before \(t(q_i,\omega)\) is the passage time of the edge \(q_i\) for the realization \(\omega.\) For any fixed path~\(\pi = (e_1,\ldots,e_r)\) having the origin as an endvertex and containing \(r\) edges, we define the truncated passage time
\begin{equation}\label{t_pi_p_def}
\hat{T}^{(k)}(\pi) = \hat{T}^{(k)}(\pi,\omega)  = \sum_{i=1}^{r} t^{(k)}(e_i,\omega).
\end{equation}
As before, for a fixed integer \(n \geq 1,\) we let \({\cal P}_n\) denote the set of all finite paths with endvertices origin and \((n,0,\ldots,0)\) and define
\begin{equation}\label{t_0n_def_p}
\hat{T}^{(k)}_n = \hat{T}^{(k)}_n(\omega)  = \inf_{\pi \in {\cal P}_n} \hat{T}^{(k)}(\pi,\omega)
\end{equation}
to be the \emph{minimum truncated passage time} between the origin and the point\\\((n,0,\ldots,0).\) We use the hat notation to emphasize the approximation of~\(T_n.\) As before we suppress the dependence on \(\omega\) unless specifically mentioned.

For integer \(i \geq 1,\) let \(f_i\) denotes the edge between \((i-1,0,\ldots,0)\) and \((i,0,\ldots,0).\) We collect the properties of the truncated passage times needed for future use.
\begin{Lemma}\label{sec_mom} Fix integers \(n,k \geq 1.\)
For any fixed integer \(k_1 \geq k\) we have
\begin{equation}\label{tk1k2}
\hat{T}_{n}^{(k)} \leq \hat{T}^{(k_1)}_{n} \leq T_{n}
\end{equation}
and so
\begin{equation}\label{unif_int}
\mathbb{E}\left(\frac{\hat{T}^{(k)}_n}{n}\right)^2 \leq \mathbb{E}\left(\frac{T_n}{n}\right)^2 \leq \frac{1}{n^2}\mathbb{E}\left(\sum_{i=1}^{n}t(f_i)\right)^2 \leq C
\end{equation}
for some constant \(C > 0\) not depending on \(n\) or \(k.\)
Also for any integer \(n_1 \geq n\) we have
\begin{equation}\label{sub_add}
|\hat{T}^{(k)}_{n}  - \hat{T}^{(k)}_{n_1}| \leq k^{\alpha}(n_1 - n).
\end{equation}
\end{Lemma}
The relations (\ref{tk1k2}) and (\ref{sub_add}) hold for every \(\omega \in \Omega.\)
In particular the result~(\ref{unif_int}) also implies that \(\frac{T_n}{n}\) is uniformly integrable.

%//CHAGNE WITH NEW NOTATION,,,, PRKVMM +eTC

\emph{Proof of Lemma~\ref{sec_mom}}:
To see (\ref{tk1k2}), we first note that for any edge \(h\) we have that the truncated passage times  satisfy
\[t^{(k)}(h)  = \min(t(h),k^{\alpha}) \leq t^{(k_1)}(h) \leq t(h).\] Thus for any fixed finite path \(\pi\) with endvertices origin and \((n,0,\ldots,0)\) we have that
\[\sum_{h \in \pi} t^{(k)}(h) \leq \sum_{h \in \pi} t^{(k_1)}(h)  \leq \sum_{h \in \pi} t(h).\]
In other words, the passage times of the paths (see (\ref{t_pi_p_def})) satisfy \[\hat{T}^{(k)}(\pi)  \leq \hat{T}^{(k_1)}(\pi) \leq T(\pi).\] Taking infimum over all such paths~\(\pi\) gives (\ref{tk1k2}).

To prove (\ref{unif_int}), we use the following relation.
\begin{equation}\label{sto_dom}
\hat{T}^{(k_1)}_n \leq T_n \leq \sum_{i=1}^{n}t(f_i)
\end{equation}
The first estimate  in (\ref{sto_dom}) is true from (\ref{tk1k2}) and the second inequality follows from (\ref{t_n_ineq123}) in Proposition~\ref{geo_prop}.
It therefore suffices to prove the final estimate in~(\ref{unif_int}). We use the estimate \((\sum_{i=1}^{l}a_i)^2 \leq l\sum_{i=1}^{l} a_i^2\) for positive~\(\{a_i\}\) with \(l = n\) and \(a_i = t(f_i),1 \leq i \leq n\) to get that \[\mathbb{E}\left(\sum_{i=1}^{n} t(f_i)\right)^2 \leq n \sum_{i=1}^{n}\mathbb{E}t^2(f_i) \leq C n^2\] where the final estimate follows from the moment condition~\((ii)\) of Section~\ref{intro}.

In what follows we prove (\ref{sub_add}). Let \(\hat{T}^{(k)}_{n,n_1}\) denote the minimum passage time between the vertices \((n,0,\ldots,0)\) and \((n_1,0,\ldots,0)\) defined analogously as \(\hat{T}^{(k)}_n\) for each \(k,n\) and \(n_1.\)
We first prove that
\begin{equation}\label{sub_add2}
\hat{T}^{(k)}_{n}  \leq \hat{T}^{(k)}_{n, n_1} + \hat{T}^{(k)}_{n_1} \text{ and } \hat{T}^{(k)}_{n_1} \leq \hat{T}^{(k)}_{n} + \hat{T}^{(k)}_{n,n_1}.
\end{equation}
We prove the first estimate in (\ref{sub_add2}). The proof of the other is analogous.

Let \(\pi\) be any finite path with endvertices as origin and \((n,0,\ldots,0)\) and let \(\pi'\) be any finite path with endvertices \((n,0,\ldots,0)\) and \((n_1,0,\ldots,0).\) The union \(\pi \cup \pi'\) contains a path with endvertices as origin and \((n_1,0,\ldots,0).\) Therefore by definition, we have that the minimum passage time between the origin and \((n_1,0,\ldots,0)\) satisfies
\[\hat{T}^{(k)}_{n_1} \leq \hat{T}^{(k)}(\pi) + \hat{T}^{(k)}(\pi').\] Here \(\hat{T}^{(k)}(.)\) is the truncated passage time as defined in (\ref{t_pi_p_def}). Taking infimum over all paths \(\pi\) gives
\[\hat{T}^{(k)}_{n_1} \leq \hat{T}^{(k)}_{n} + \hat{T}^{(k)}(\pi').\] Taking infimum over all paths \(\pi'\) then gives the first estimate in~(\ref{sub_add2}).

Thus from (\ref{sub_add2}) we have
\begin{equation}\label{t_est_temp123}
|\hat{T}^{(k)}_{n}  - \hat{T}^{(k)}_{n_1}| \leq \hat{T}^{(k)}_{n,n_1} \leq \sum_{i=n+1}^{n_1} t^{(k)}(f_i)
\end{equation}
where as before \(f_i\) denotes the edge from \((i-1,0,\ldots,0)\) to \((i,0,\ldots,0).\) The final estimate again follows from the definition of the minimum passage time. Since for any edge \(h,\) the truncated passage time satisfies \[t^{(k)}(h) = \min(t(h),k^{\alpha}) \leq k^{\alpha},\] we have from the last estimate in (\ref{t_est_temp123}) that
\[|\hat{T}^{(k)}_{n}  - \hat{T}^{(k)}_{n_1}| \leq k^{\alpha}(n_1-n).\] This proves (\ref{sub_add}).~\(\qed\)

%We let \(\hat{T}^{(n)}_k\) be the shortest passage time from \((0,...,0)\) to \((k,0,...,0)\) when the passage times are \(t_n(e_i).\) %= \min(t(e_i), n^{\alpha}).\)

%Since \(t^{(k)}(e) \leq t(e)\) for any edge \(e,\) we have that the moment conditions (i) and (ii) are also satisfied by \(\{t^{(k)}(q_i)\}_{i \geq 1}\) for any fixed integer \(k \geq 1.\)

\subsection*{Geodesics}
To study further the relation between the truncated and untruncated passage times, we need the concept of geodesics. In the first part of this subsection, we study the existence of geodesics and the second part we derive auxiliary estimates regarding the geodesics needed for future use.

We have the following main result regarding the existence of geodesics. As before, for integer \(i \geq 1,\) we let \(f_i\) denote the edge between~\((i-1,0,\ldots,0)\) and~\((i,0,\ldots,0).\)
\begin{Proposition}\label{geo_prop} There is a null set \(Z_0 \subseteq \Omega\) with \(\mathbb{P}(Z_0) = 0\) so the following holds for all \(\omega \in Z_0^c.\) For all integers \(n \geq 1,\) there are unique (finite) paths \(\pi_n = \pi_n(\omega)\) and \(\hat{\pi}^{(n)}_n = \hat{\pi}^{(n)}_n(\omega)\) such that
\begin{equation}\label{geo_orig}
T_n = T(\pi_n)
\end{equation}
and
\begin{equation}\label{geo_trunc}
\hat{T}^{(n)}_n = \hat{T}^{(n)}(\hat{\pi}^{(n)}_n).
\end{equation}
If \(\omega \in Z_0,\) then for every integer \(n \geq 1,\) we define \(\pi_n = \hat{\pi}^{(n)}_n = (f_1,\ldots,f_{n}),\) i.e., the straight line with endvertices as origin and the point \((n,\ldots,0).\)\\
For all \(\omega \in \Omega,\) we therefore have
\begin{equation}\label{t_n_ineq123}
T_n \leq T(\pi_n) \leq \sum_{i=1}^{n} t(f_i) \text{  and  }   \hat{T}^{(n)}_n \leq \hat{T}^{(n)}(\hat{\pi}^{(n)}_n) \leq \sum_{i=1}^{n}t^{(n)}(f_i).
\end{equation}
\end{Proposition}
We refer to \(\pi_n\) and \(\hat{\pi}^{(n)}_n\) as geodesics. In what follows, we also derive certain auxiliary estimates along the way that are useful for estimating the variances in Section~\ref{pf1}.

The first step is to see that long paths have sizeable passage time with high probability. Let
\begin{equation}\label{mu_def}
\mu := \sup_{i}\mathbb{E}t(q_i).
\end{equation}
Using conditions \((i)\) and \((ii)\) in Section~\ref{intro}, we have that
\begin{equation}
\mu \in (0,\infty).
\end{equation}
To see that \(\mu < \infty,\) we use Cauchy-Schwarz inequality and have for any edge~\(e\) that \[\mathbb{E}t(e) \leq \left(\mathbb{E}t^2(e)\right)^{1/2}.\] Using condition \((ii)\) in Section~\ref{intro}, we therefore have \(\mu < \infty.\) For \(\mu > 0,\) we fix edge \(e\) and \(\epsilon > 0\) and have that \[\mathbb{E}t(e) \geq \epsilon \mathbb{P}(t(e) \geq  \epsilon).\] Using condition \((i)\) of Section~\ref{intro}, we therefore have that \[\mu \geq \inf_i \mathbb{E}t(q_i) \geq \epsilon \inf_i\mathbb{P}(t(q_i) \geq \epsilon) \geq \frac{\epsilon}{2},\] if \(\epsilon > 0\) is small.

For any fixed path \(\pi = (e_1,\ldots,e_m)\) containing \(m\) edges, let~\(T(\pi)\) and~\(\hat{T}^{(k)}(\pi)\) be the untruncated and truncated passage times as defined in~(\ref{t_pi_def}) and~(\ref{t_pi_p_def}), respectively.
We have the following Lemma.
\begin{Lemma}\label{t_pi_est} Let \(\mu \in (0,\infty)\) be as in (\ref{mu_def}). We have
\begin{equation}\label{t_pi}
\mathbb{P}(T(\pi) \leq \beta_1 m) \leq e^{-dm}
\end{equation}
and
\begin{equation}\label{t_pi_p}
\mathbb{P}( \hat{T}^{(k)}(\pi) \leq \beta_1 m) \leq e^{-dm}
\end{equation}
for some positive constant \(\beta_1 \in (0, \mu)\) and for all \(m, k \geq 1.\)
\end{Lemma}
Here and henceforth all constants are independent of \(n\) and \(k.\)\\
\emph{Proof Lemma~\ref{t_pi_est}}: To prove (\ref{t_pi}), we write \[\mathbb{P}(T(\pi) \leq \beta m)  = \mathbb{P}\left(\sum_{i=1}^{m}t(e_i) \leq \beta m\right)\] for a fixed \(\beta > 0.\) For a fixed \(s  >0\) we have that \(\sum_{i=1}^{m} t(e_i) \leq \beta m\) if and only if the product \(\prod_{i=1}^{m} e^{-st(e_i)} \geq e^{-s\beta m}\) and so
\begin{eqnarray}
\mathbb{P}(T(\pi) \leq \beta m) &=& \mathbb{P}\left(\prod_{i=1}^{m} e^{-st(e_i)} \geq e^{-s\beta m}\right)  \nonumber\\
&\leq& e^{s\beta m} \mathbb{E}\left(\prod_{i=1}^{m}e^{-st(e_i)}\right) \nonumber\\
&=& e^{s\beta m} \prod_{i=1}^{m}\mathbb{E}\left(e^{-st(e_i)}\right) \label{y_1_eq1}
\end{eqnarray}
where the first inequality follows from the Markov inequality and the equality~(\ref{y_1_eq1}) follows since \(\{t(e_{i})\}_{1 \leq i \leq m}\) are independent.

For a fixed \(\epsilon > 0,\) we have that
\begin{eqnarray}
\mathbb{E}e^{-st(e_i)} &=& \int_{t(e_i) < \epsilon} e^{-st(e_i)} d\mathbb{P} + \int_{t(e_i) \geq \epsilon} e^{-st(e_i)} d\mathbb{P} \nonumber\\
&\leq& \int_{t(e_i) < \epsilon} e^{-st(e_i)} d\mathbb{P} + e^{-s\epsilon} \label{t_pi_an}\\
&\leq& \mathbb{P}(t(e_i) < \epsilon) + e^{-s\epsilon}. \nonumber
\end{eqnarray}
Thus for any fixed \(i \geq 1,\) we have
\[\mathbb{E}e^{-st(e_i)} \leq \sup_{j \geq 1} \mathbb{E}e^{-st(e_j)} \leq \sup_j \mathbb{P}(t(e_j) < \epsilon) + e^{-s\epsilon}.\]
Using condition \((i)\) of Section~\ref{intro}, the first term in the last expression is less than \(\frac{e^{-6d}}{2}\) if \(\epsilon  = \epsilon(d) > 0\) is small. Fixing such an \(\epsilon,\) we choose~\(s = s(\epsilon) > 0\) large so that the second term is also less than~\(\frac{e^{-6d}}{2}.\) The choices of~\(\epsilon\) and~\(s\) do not depend on the index~\(i.\) Therefore substituting into (\ref{y_1_eq1}), we have that \[\mathbb{P}(T(\pi) \leq \beta m) \leq e^{s\beta m} e^{-6dm} \leq e^{-2dm},\] for all \(m \geq 1\) provided \(\beta = \beta(s,d) > 0\) is small enough so that \(s\beta < 4d.\) Without loss of generality, we let \(\beta < \mu.\)

The proof of (\ref{t_pi_p}) is analogous. We use the fact that \(\{t^{(k)}(e_{i})\}_{1 \leq i \leq m}\) are independent and thus for a fixed \(s  >0\) we have (as in (\ref{y_1_eq1})) that
\begin{equation}\label{y_2_eq2}
\mathbb{P}(\hat{T}^{(k)}(\pi) \leq \beta m) = \mathbb{P}\left(\sum_{i=1}^{m} t^{(k)}(e_i) \leq \beta m\right) \leq e^{s\beta m} \prod_{i=1}^{m}\mathbb{E}(e^{-st^{(k)}(e_i)}).
\end{equation}
For a fixed \(0 < \epsilon <1,\) we have that
\begin{eqnarray}
\mathbb{E}e^{-st^{(k)}(e_i)} &=& \int_{t^{(k)}(e_i) < \epsilon} e^{-st^{(k)}(e_i)} d\mathbb{P} + \int_{t^{(k)}(e_i) \geq \epsilon} e^{-st^{(k)}(e_i)} d\mathbb{P} \nonumber\\
&\leq& \int_{t^{(k)}(e_i) < \epsilon} e^{-st^{(k)}(e_i)} d\mathbb{P} + e^{-s\epsilon} \nonumber\\
&=& \int_{t(e_i) < \epsilon} e^{-st(e_i)} d\mathbb{P} + e^{-s\epsilon} \nonumber
\end{eqnarray}
which is the same as (\ref{t_pi_an}). The final equality is because \(\epsilon < 1\) and thus \(t^{(k)}(e_i)< \epsilon\) if and only if \(t(e_i) < \epsilon.\) By an analogous analysis following (\ref{t_pi_an}) we obtain (\ref{t_pi_p}). \(\qed\)

The next step in the construction is to obtain estimates on the straight line joining origin to \((n,0,\ldots,0).\) As in the statement of the Proposition~\ref{geo_prop}, for \(i \geq 1,\) let \(f_i\) denote the edge between \((i-1,0,\ldots,0)\) and \((i,0,\ldots,0)\) and let
\begin{equation}\label{a_n_def}
A_n = \left\{\sum_{i=1}^{n} t(f_i) \leq 2\mu n\right\},
\end{equation}
and for \(k \geq 1,\) let
\begin{equation}\label{a_n_p_def}
{A}_n(k) = \left\{\sum_{i=1}^{n} t^{(k)}(f_i) \leq 2\mu n\right\},
\end{equation}
where as before \(\mu \in (0,\infty)\) is as in (\ref{mu_def}).
We have the following result.
\begin{Lemma}\label{lemma_a_n_est2} We have
\begin{equation}\label{an_prob2}
\mathbb{P}\left(\liminf_n A_n\right) = 1
\end{equation}
and
\begin{equation}\label{an_prob_p2}
\mathbb{P}\left(\liminf_n {A}_n(n)\right) =1.
\end{equation}
\end{Lemma}
\emph{Proof of Lemma~\ref{lemma_a_n_est2}}: We prove (\ref{an_prob2}) first. Letting \(X_i = t(f_i) - \mathbb{E}t(f_i),\) we have that
\begin{equation}\label{sum_var}
\sum_{i \geq 1} \frac{var(X_i)}{i^2} \leq \sum_{i \geq 1} \frac{\mathbb{E}t^2(f_i)}{i^2} \leq C_1 \sum_{i \geq 1}\frac{1}{i^2} < \infty
\end{equation}
where \(C_1 > 0\) is a constant and the second inequality follows from the moment condition \((ii)\) of Section~\ref{intro}.
Therefore using Kronecker's Lemma (Theorem \(2.5.5,\) Chapter \(2,\) Durrett~(2010)), we have that
\begin{equation}\label{x_i_conv}
\frac{1}{n}\sum_{i=1}^{n}X_i = \frac{1}{n}\sum_{i=1}^{n}t(f_i) -\frac{1}{n}\sum_{i=1}^{n}\mathbb{E}t(f_i) \longrightarrow 0 \text{ a.s. }
\end{equation}
as \(n \rightarrow \infty.\) Let \(Z\) denote the null set in (\ref{x_i_conv}) with \(\mathbb{P}(Z) = 0\) and fix \(\omega \in Z^c.\) There exists \(N = N(\omega)\) so that for all \(n \geq N\) we have
\begin{equation}\label{x_i_conv2}
\frac{1}{n}\sum_{i=1}^{n}t(f_i) \leq \frac{1}{n}\sum_{i=1}^{n}\mathbb{E}t(f_i)  + \mu \leq 2\mu.
\end{equation}
The final estimate holds since \(\mathbb{E}t(f_i) \leq \mu\) by definition. This implies that \(\omega \in A_n\) for all \(n \geq N.\) Thus \(\omega \in \liminf_n A_n\) and so (\ref{an_prob2}) holds.

For the other estimate (\ref{an_prob_p2}), we have that \(t^{(n)}(f_i) \leq t(f_i)\) for all \(i\) and therefore \(A_n \subseteq A_n(n)\) for all \(n\) and so (\ref{an_prob_p2}) also holds. \(\qed\)

As a final step before proving Proposition~\ref{geo_prop}, we need lower bounds on passage times of long paths. Let \(\beta_1 > 0\) be as in Lemma~\ref{t_pi_est} and let
\begin{equation}\label{e_k_def}
E_{m} := \bigcup_{r \geq \frac{8\mu}{\beta_1} m}\;\;\bigcup_{\pi \in {\cal Q}_r} \left\{T(\pi) < \beta_1 r\right\}
\end{equation}
denote the event that there exists a path \(\pi\) with origin as an endvertex and consisting of \(r \geq \frac{8 \mu}{\beta_1}m\) edges, whose passage time~\(T(\pi)\) defined in (\ref{t_pi_def}) is less than~\(\beta_1 r.\) In the above \({\cal Q}_r\) denotes the set of all paths with origin as an endvertex and consisting of~\(r\) edges. Since there are at most \((2d)^{r}\) paths in~\({\cal Q}_r,\) we have using (\ref{t_pi}) that
\begin{equation}\label{a_0k}
\mathbb{P}(E_{m}) \leq \sum_{r \geq 8\mu\beta_1^{-1} m} (2d)^{r}e^{-dr} \leq \sum_{r \geq 8\mu\beta_1^{-1} m} e^{-\beta_{22} r} \leq C_2 e^{-\beta_2 m}
\end{equation}
for all \(m \geq 1\) and for some positive constants \(\beta_2, \beta_{22}\) and \(C_2.\) Here we use \(xe^{-x} \) attains its maximum at \(x = 1\) and so \(2d e^{-d} \leq 2e^{-1} =: e^{-\beta_{22}}\) for all \(d \geq 2.\)

We define a similar event for the truncated random variables. For \(m \geq 1,\) let
\begin{equation}\label{e_k_hat_def}
{E}_m(k) := \bigcup_{r \geq \frac{8\mu}{\beta_1} m}\;\;\bigcup_{\pi \in {\cal Q}_r} \left\{\hat{T}^{(k)}(\pi) < \beta_1 r\right\}
\end{equation}
denote the event that there exists a path \(\pi\) starting from the origin containing \(r \geq \frac{8 \mu}{\beta_1}m\) edges and whose truncated passage time \(\hat{T}^{(k)}(\pi)\) defined in (\ref{t_pi_p_def}) is less than~\(\beta_1 r.\) Using (\ref{t_pi_p}) and proceeding as in (\ref{a_0k}), we have that
\begin{equation}\label{a_0k_p}
\mathbb{P}({E}_{m}(k)) \leq C_2 e^{-\beta_2 m}
\end{equation}
for all \(m, k \geq 1,\)  where \(C_2, \beta_2 > 0\) are the constants in (\ref{a_0k}).

Using (\ref{a_0k}) and (\ref{a_0k_p}) individually with the Borel-Cantelli lemma, we therefore have the following result.
%\begin{Lemma}\label{lemma_en_est}
We have
\begin{equation}\label{en_prob2}
\mathbb{P}\left(\liminf_n E^c_n\right) = 1
\end{equation}
and
\begin{equation}\label{en_prob_p2}
\mathbb{P}\left(\liminf_n {E}^c_n(n)\right) =1.
\end{equation}
%\end{Lemma}

We collect the above results to prove Proposition~\ref{geo_prop}. For \(m \geq 1,\) let
\begin{equation}\label{b_m_def}
B_m := [-m,m]^d
\end{equation}
denote the box with side length \(2m.\)\\
\emph{Proof of Proposition~\ref{geo_prop}}:  We prove (\ref{geo_orig}) and an analogous proof holds for~(\ref{geo_trunc}). Let \(E_n\) and \(A_n\) be the events defined in (\ref{e_k_def}) and (\ref{a_n_def}), respectively.
Setting
\begin{equation}\label{f_n_def}
F_n = E_{n}^c \cap A_n,
\end{equation}
we have from (\ref{an_prob2}) and (\ref{en_prob2}) that
\begin{equation}\label{f_n_est}
\mathbb{P}\left(\liminf_n F_n \right) = 1.
\end{equation}
Suppose \(\omega \in \liminf_n F_n.\) There exists \(N_1 = N_1(\omega)\) such that for all \(n \geq N_1,\) we have \(\omega \in F_n = E_n^c \cap A_n.\) Fix \(n \geq N_1.\) Since \(\omega \in A_n,\) the minimum passage time \(T_n = T_n(\omega)\) between the origin and the point~\((n,0,\ldots,0)\) is at most~\(2\mu n.\) Since \(\omega \in E^c_{n},\) every path containing the origin as an endvertex and consisting of \(r \geq \frac{8\mu}{\beta_1}n\) edges, has untruncated passage time of at least~\(\beta_1 r \geq 8\mu n.\) Therefore, some path~\(\pi\) contained completely in \(B_{8\mu\beta_1^{-1}n}\) has passage time \(T(\pi)   = \sum_{e \in \pi} t(e) \leq 2\mu n.\) In particular, there exists a path \(\pi_n = \pi_n(\omega)\) contained in \(B_{8\mu\beta_1^{-1}n}\) such that
\begin{equation}\label{t_i_est}
T_n = T(\pi_n).
\end{equation}
If there are multiple choices for any \(\pi_n,\) we use a fixed rule, for example, the iterative procedure described at the end of this section to choose a path.

If \(n < N_1,\) again we have \(\sum_{i=1}^{n}t(f_i) \leq \sum_{i=1}^{N_1}t(f_i) \leq 2\mu N_1,\) since \(\omega~\in~A_{N_1}.\) Also, since \(\omega \in E^c_{N_1},\) every path containing the origin as an endvertex and consisting of \(r \geq \frac{8\mu}{\beta_1}N_1\) edges has passage time at least \(\beta_1 r \geq 8\mu N_1.\) Arguing as before, some path~\(\pi\) contained completely in~\(B_{8\mu\beta_1^{-1}N_1}\) has passage time~\(T(\pi)\) less than or equal to~\(2\mu N_1\) and so there exists a path \(\pi_n\) satisfying~(\ref{t_i_est}) and completely contained in~\(B_{8\mu\beta_1^{-1}N_1}.\)

If \(\omega \notin \liminf_n F_n,\) then \(\omega\) belongs to a null set. We set \(\pi_n\) to be the path consisting of the edges \(\{f_i\}_{1 \leq i \leq n-1}\) and therefore \(\pi_n\) is the straight line joining the origin to \((n,0,\ldots,0).\) By construction, we have that (\ref{t_n_ineq123}) is satisfied.

We do an analogous analysis as above for the truncated random variable~\(\hat{T}^{(n)}_n.\)
Recalling the events \(E_n(n)\) and \(A_n(n)\) defined in (\ref{e_k_hat_def}) and (\ref{a_n_p_def}), respectively, we set
\begin{equation}\label{f_n_hat_def}
\hat{F}_n = E_{n}^c(n) \cap A_n(n),
\end{equation}
we have from (\ref{an_prob_p2}) and (\ref{en_prob_p2}) that
\begin{equation}\label{f_n_est_hat}
\mathbb{P}\left(\liminf_n \hat{F}_n \right) = 1.
\end{equation}
We then obtain the geodesic \(\hat{\pi}^{(n)}_n\) as in the discussion in the paragraph following (\ref{f_n_est}).

Finally, we set \[Z_0^c := \liminf_n (F_n \cap \hat{F}_n)\] and obtain from (\ref{f_n_est}) and (\ref{f_n_est_hat}) that \(\mathbb{P}(Z_0) = 0.\) And for all \(\omega \in Z_0^c,\) finite paths attain the minimum truncated and untruncated passage times as described above.

We prove the first set of inequalities in~(\ref{t_n_ineq123}). An analogous proof holds for the other set. Let \(\pi_0\) denotes the straight line with origin and \((n,0,\ldots,0)\) as endvertices. Fix \(\omega \in Z_0^c.\) We then have \(T(\pi_n) = T_n \leq T(\pi_0)\) where the final inequality is by definition of the minimum passage time (see (\ref{t_0n_def})). If \(\omega \in Z_0,\) then \(T_n \leq T(\pi_n) = T(\pi_0)\) where the first inequality again is by the definition of minimum passage time in (\ref{t_0n_def}). \(\qed\)

%as follows: for \(k = K_0(\omega),\) define \(\hat{\pi}^{(n)}_k(\omega)\) as in (\ref{t_m_i_est}) for \(1 \leq n \leq 2k.\) For \(k = K_0(\omega)+1,\) define  \hat{\pi}^{(n)}_k(\omega)\) as in (\ref{t_m_i_est}) for \(1 \leq n \leq 2k.\) %If there is more than one path that attains the shortest time, we provide an iterative procedure at the end of this section to choose a unique path. Similarly, for every \(n \geq 1\) and for every \(1 \leq k \leq 2n,\) define \(\hat{T}^{(n)}_k = \hat{T}^{(n)}_k(\omega)\) to be the shortest time taken for reaching \((k,0)\) from \((0,0).\)

 %An analogous estimate holds for \(V_2\) and this completes the proof of (\ref{t_pi}).
Finally, for completeness, we provide an iterative procedure to choose a single path in the presence of multiple choices. We remark that it is also possible to choose a path using any deterministic rule like for e.g., a fixed ordering of the paths. For simplicity we consider the case for \(d~=~2.\) An analogous procedure holds for general \(d.\) Let \({\cal S}_1 = \{L_i\}_{1 \leq i \leq W} =\{(S_{i,1},...,S_{i,H_i})\}_{1 \leq i \leq W}\) be any finite set of paths with endvertices \((0,0)\) and \((n,0).\)

Let \(x_{i,j}\) and \(y_{i,j}\) be the \(x\)- and \(y\)-coordinates, respectively, of the centre of the edge \(S_{i,j}.\) Let \(y'_1 = \min_{L_k \in {\cal S}_1} y_{k,1}\) and let \({\cal S}'_1 = \{L_k \in {\cal S}_1 : y_{k,1} = y'_1\}.\) Let \(x'_1 = \min_{L_k \in {\cal S}'_1} x_{k,1}.\) Let \(h_1\) be the edge attached to the origin whose centre has coordinates \((x'_1,y'_1).\) Clearly \(h_1\) is the first edge of some path in~\({\cal S}'_1.\) Let \({\cal S}_2\) be the set of paths in \({\cal S}'_1\) whose first edge is \(h_1.\) Repeating the above procedure with \({\cal S}_2,\) we obtain an edge \(h_2\) attached to \(h_1.\) Continuing iteratively, this procedure terminates after a finite number of steps resulting in a unique path. Also, the final path obtained does not depend on the initial ordering of the paths.

\subsubsection*{Geodesics contained in large boxes}
The following result estimates the probability that the geodesics \(\pi_n\) and \(\hat{\pi}^{(n)}_n\) are contained in large finite boxes.
\begin{Lemma}\label{box_est} Let \(\mu = \sup_i \mathbb{E}t(e_i) \in (0,\infty)\) be as in~(\ref{mu_def}) and let \(\beta_1 \in (0, \mu)\) be the constant in Lemma~\ref{t_pi_est}. Fix integer \(n \geq 1\) and \(\epsilon > 0\) a constant and let \(n_{\epsilon} = n^{1+\epsilon}.\) We have
\begin{equation}\label{box_pi}
\mathbb{P}\left(\pi_n \subseteq B_{8\mu \beta_1^{-1} n_{\epsilon}}\right) \geq 1 - \frac{C}{n^{1+2\epsilon}}
\end{equation}
and
\begin{equation}\label{box_pi_hat}
\mathbb{P}\left(\hat{\pi}^{(n)}_n \subseteq B_{8\mu \beta_1^{-1} n_{\epsilon}}\right) \geq 1 - \frac{C}{n^{1+2\epsilon}}
\end{equation}
for some constant \(C > 0\) and for all \(n \geq 1.\)
\end{Lemma}
In the above we use the notation \(\left\{\pi_n \subseteq B_{8\mu \beta_1^{-1} n_{\epsilon}}\right\}\) to denote the event that every edge of \(\pi_n\) is contained in the box \(B_{8\mu \beta_1^{-1} n_{\epsilon}}.\)\\

\emph{Proof of Lemma~\ref{box_est}}: We prove (\ref{box_pi_hat}) and the proof is analogous for (\ref{box_pi}). As in Proposition~\ref{geo_prop}, let \(f_i\) denote the edge between \((i-1,0,\ldots,0)\) and \((i,0,\ldots,0).\) Let \(E_n(k)\) be the event defined in (\ref{e_k_hat_def}) and similar to (\ref{a_n_p_def})
define the event
\begin{equation}\label{a_npp_def}
\hat{A}_{n_{\epsilon}}(n) = \left\{\sum_{i=1}^{n}t^{(n)}(f_i) \leq 2\mu n^{1+\epsilon} \right\}
\end{equation}
Setting
\begin{equation}\label{h_n_def}
H_n = E_{n_{\epsilon}}^c(n) \cap \hat{A}_{n_{\epsilon}}(n),
\end{equation}
we suppose that \(H_n\) occurs. The minimum truncated passage time \(\hat{T}^{(n)}_n\) between the origin and the point~\((n,0,\ldots,0)\) is at most \(2\mu n^{1+\epsilon}.\) Since \(E_{n_{\epsilon}}^c(n)\) also occurs, every path \(\pi\) starting from the origin and containing \(r \geq \frac{8\mu}{\beta_1}n^{1+\epsilon}\) edges has truncated passage time \(\hat{T}^{(n)}(\pi) \geq \beta_1 r \geq 8\mu n^{1+\epsilon}.\) Therefore arguing as in the paragraph preceding (\ref{t_i_est}), we have that the geodesic~\(\hat{\pi}^{(n)}_n\) with endvertices as origin and \((n,0,\ldots,0)\) is contained in~\(B_{8\mu\beta^{-1}_1n_{\epsilon}}.\)

To estimate \(\mathbb{P}(\hat{A}_{n_{\epsilon}}(n)),\) we first have
\[\mathbb{P}(\hat{A}^c_{n_{\epsilon}}(n))  = \mathbb{P}\left(\sum_{i=1}^{n}t^{(n)}(f_i)  > 2\mu n^{1+\epsilon}\right) \leq \mathbb{P}\left(\sum_{i=1}^{n}t(f_i)  > 2\mu n^{1+\epsilon}\right)\] since \(t^{(n)}(f_i) \leq t(f_i)\) for all \(i\) (see (\ref{trunc_pas})).
Letting \(X_i  = t(f_i) - \mathbb{E}t(f_i)\) and using \(\mathbb{E}t(f_i) \leq \mu\) for all \(i,\) we have that
\[\mathbb{P}(\hat{A}^c_{n_{\epsilon}}(n)) \leq \mathbb{P}\left(\sum_{i=1}^{n}X_i \geq 2\mu n^{1+\epsilon} - \mu n\right)\leq \mathbb{P}\left(\sum_{i=1}^{n}X_i \geq \mu n^{1+\epsilon}\right).\]
Using Markov inequality on the last estimate, we then have
\begin{equation}\label{hat_an_est}
\mathbb{P}(\hat{A}^c_{n_{\epsilon}}(n)) \leq \frac{\mathbb{E}\left(\sum_{i=1}^{n}X_i\right)^2}{\mu^2 n^{2+2\epsilon}} \leq \frac{C_1 n}{n^{2+2\epsilon}}
\end{equation}
for some constant \(C_1 > 0,\) where the final estimate follows from
\begin{eqnarray}\label{x_i_est2}
\mathbb{E}\left(\sum_{i=1}^{n} X_{i}\right)^{2} = \sum_{i=1}^{n}\mathbb{E}X_{i}^2 + \sum_{i \neq j} \mathbb{E}X_{i}X_{j} =\sum_{i=1}^{n} \mathbb{E}X_i^2 \leq C_1n
\end{eqnarray}
for some constant \(C_1 > 0.\) The second equality in (\ref{x_i_est2}) follows from the independence of \(X_i\) and \(X_j\) for \(j \neq i\) so that \(\mathbb{E}X_iX_j = \mathbb{E}X_i\mathbb{E}X_j = 0.\) The final estimate in (\ref{x_i_est2}) follows from the moment condition \((ii)\) in Section~\ref{intro}.

From (\ref{a_0k}) and (\ref{hat_an_est}), we have that
\begin{equation} \label{fn_prob}
\mathbb{P}(H_n^c) \leq \mathbb{P}(E_{n_{\epsilon}}(n)) + \mathbb{P}(\hat{A}^c_{n_{\epsilon}}(n)) \leq C_2 e^{-\beta_2 n^{1+\epsilon}} + \frac{C_1}{n^{1+2\epsilon}} \leq \frac{C_3}{n^{1+2\epsilon}}
\end{equation}
for some constant \(C_3 > 0.\) This proves (\ref{box_pi_hat}) and an analogous proof holds for (\ref{box_pi}).~\(\qed\)

%&=& n-1 + \frac{1}{1-e^{-\beta_2}}e^{-\beta_2 n}  \nonumber\\
%&\leq& n \nonumber

%and
%\begin{equation}\label{w_n_est_tot}
%\mathbb{P}(W_n^c) \leq \frac{C}{n^{2+2\delta}}
%\end{equation}

\setcounter{equation}{0}
\renewcommand\theequation{\thesection.\arabic{equation}}
\section{Truncated and untruncated passage times are asymptotically equivalent}\label{pf1}
Fix integer \(n \geq 1.\) To see that the truncated minimum passage time~\(\hat{T}^{(n)}_n\) is a good approximation of \(T_n,\) we need the estimates on the probability that \(T_n = \hat{T}^{(n)}_n.\) Defining the events
\begin{equation}\label{v_n_def}
V_n = \{T_n = \hat{T}^{(n)}_n\} \text{ and } W_n = \{\hat{T}^{((n+1)^2)}_{n^2} = \hat{T}^{(n^2)}_{n^2}\},
\end{equation}
we have the following result.
\begin{Proposition} \label{lemma_v_n} Let \(\delta = \frac{1}{2d}\) where \(d \geq 2\) is the dimension. We have
\begin{equation}\label{v_n_est_tot}
\mathbb{P}(V_n^c) \leq \frac{C}{n^{1+\delta}}
\end{equation}
and
\begin{equation}\label{w_n_est_tot}
\mathbb{P}(W_n^c) \leq \frac{C}{n^{2+2\delta}}
\end{equation}
for some constant \(C > 0\) and for all \(n \geq 1.\) In particular, we have
\begin{equation}\label{t_n_eq_hat}
\frac{1}{n}\left(T_n - \mathbb{E}T_n\right) - \frac{1}{n}\left(\hat{T}^{(n)}_n - \mathbb{E}\hat{T}^{(n)}_n\right) \longrightarrow 0\text{ a.s.}
\end{equation}
as \(n \rightarrow \infty.\)
\end{Proposition}
The above result implies that it suffices to study the convergence of \(\frac{1}{n}\left(\hat{T}^{(n)}_n - \mathbb{E}\hat{T}^{(n)}_n\right).\)

In what follows we first derive preliminary results related to \(V_n\) and \(W_n\) and then finally prove Proposition~\ref{lemma_v_n}. Let~\(Z_0\) be the null set in Proposition~\ref{geo_prop} so that for all \(\omega \in Z_0^c,\) the minimum passage times (truncated and untruncated) are attained by finite paths. Fix \(\omega \in Z_0^c\) and integer \(n \geq 1\) and let \(\hat{\pi}^{(n)}_n = \hat{\pi}^{(n)}_n(\omega)\) be the geodesic for \(\hat{T}^{(n)}_n\) as in Proposition~\ref{geo_prop} so that the truncated passage time \(\hat{T}^{(n)}(.)\) as defined in~(\ref{t_pi_p_def}) satisfies \begin{equation}\label{eq_tth2}
\hat{T}^{(n)}_n = \hat{T}^{(n)}(\hat{\pi}^{(n)}_n).
\end{equation}

The geodesic \(\hat{\pi}^{(n)}_n\) has finite number of edges and the following crucial observation regarding the passage times is a first step to estimate the probability of the event~\(V_n.\)
\begin{Lemma}\label{cruc}Fix \(\omega \in Z_0^c.\) If \(t(e) < n^{\alpha}\) for all \(e \in \hat{\pi}^{(n)}_n(\omega),\) then \(T_n(\omega) = \hat{T}^{(n)}_n(\omega);\) i.e. \(\omega \in V_n.\) Here \(\alpha >0\) is as in~(\ref{trunc_pas}).
\end{Lemma}
In words, if the untruncated passage time of each edge in the geodesic~\(\hat{\pi}^{(n)}_n\) is less than~\(n^{\alpha},\) then the truncated and the untruncated minimum passage times are equal.

\emph{Proof of Lemma~\ref{cruc}}: Fix \(\omega \in Z_0^c\) and let \(L = \#\hat{\pi}^{(n)}_n\) denote the number of edges in \(\hat{\pi}^{(n)}_n = \hat{\pi}^{(n)}_n(\omega)\) and let \(\hat{\pi}^{(n)}_n = (h_1,\ldots,h_L)\) denote the edges of the geodesic~\(\hat{\pi}^{(n)}_n.\) Since \(t(h_i) < n^{\alpha}\) we have that the truncated passage time \(t^{(n)}(h_i) = \min(t(h_i),n^{\alpha}) = t(h_i)\) and so
\begin{equation}\label{eq_tth}
T(\hat{\pi}^{(n)}_n) = \sum_{i=1}^{L}t(h_i) = \sum_{i=1}^{L} t^{(n)}(h_i) = \hat{T}^{(n)}(\hat{\pi}^{(n)}_n) = \hat{T}^{(n)}_n.
\end{equation}
The first equality follows from the definition of the untruncated path passage time \(T(.)\) in (\ref{t_pi_def}). The third equality follows from the definition of the truncated passage time~\(\hat{T}^{(n)}(.)\) in~(\ref{t_pi_p_def}) and the final equality follows from~(\ref{eq_tth2}).

For any fixed path \(\pi = (g_1,\ldots,g_r)\) with endvertices as origin and the point \((n,0,\ldots,0),\) we therefore have
\begin{equation}\label{eq_tth3}
T(\hat{\pi}^{(n)}_n) = \hat{T}^{(n)}(\hat{\pi}^{(n)}_n) = \hat{T}^{(n)}_n \leq \hat{T}^{(n)}(\pi) =\sum_{i=1}^{r}t^{(n)}(g_i) \leq \sum_{i=1}^{r}t(g_i) = T(\pi).
\end{equation}
The first equality follows from (\ref{eq_tth}). The second equality follows from (\ref{eq_tth2}) since \(\omega \in Z_0^c.\) The first inequality follows from the definition of the minimum truncated passage time in~(\ref{t_0n_def_p}). The second inequality follows since \(t^{(n)}(g_i) \leq t(g_i)\) for all~\(n \geq 1\) (see (\ref{trunc_pas})).

Taking infimum over all finite paths~\(\pi\) in (\ref{eq_tth3}) and using the definition of minimum passage time in (\ref{t_0n_def}), we have \[T(\hat{\pi}^{(n)}_n) \leq T_n \leq T(\hat{\pi}^{(n)}_n)\] where the final estimate holds by the definition of the minimum passage time. But this means that equality holds in the above expression and we have \[T_n = T(\hat{\pi}^{(n)}_n) = \hat{T}^{(n)}(\hat{\pi}^{(n)}_n)  = \hat{T}^{(n)}_n.\] The middle equality follows from~(\ref{eq_tth}) above. This proves the Lemma.~\(\qed\)

The observation in Lemma~\ref{cruc} along with the fact that geodesics are contained in finite boxes with high probability (see Lemma~\ref{box_est}) allows us to estimate the event~\(V_n.\) Let \(\epsilon_0 > 0\) be a constant (to be determined later) and recalling the notation in Lemma~\ref{box_est}, let
\begin{equation}\label{g_n_def}
G_n = \left\{\hat{\pi}^{(n)}_n \subseteq B_{8\mu \beta_1^{-1} n^{1+\epsilon_0}}\right\}
\end{equation}
denote the event that every edge of the geodesic~\(\hat{\pi}^{(n)}_n \) is contained in the box~\(B_{n_0} = [-n_0,n_0]^d,\) where \(n_0 = 8\mu \beta_1^{-1} n^{1+\epsilon_0}.\) Here \(\mu = \sup_i \mathbb{E}t(e_i) \in (0,\infty)\) is as in~(\ref{mu_def}) and \(\beta_1 \in (0, \mu)\) is the constant in Lemma~\ref{t_pi_est}.

We have
\begin{eqnarray}
\mathbb{P}(V_n^c) &=& \mathbb{P}(V^c_n \cap G_n) + \mathbb{P}(V_n^c \cap G_n^c) \nonumber\\
&\leq& \mathbb{P}(V^c_n \cap G_n) + \mathbb{P}(G_n^c) \nonumber\\
&\leq& \mathbb{P}(V^c_n \cap G_n) + \frac{C_2}{n^{1+2\epsilon_0}} \label{v_n_est1}\\
&=& \mathbb{P}(V^c_n \cap G_n \cap Z_0^c) + \frac{C_2}{n^{1+2\epsilon_0}} \label{v_n_est22}
\end{eqnarray}
for some constant \(C_2 > 0,\) where (\ref{v_n_est1}) follows using (\ref{box_pi_hat}) and (\ref{v_n_est22}) follows since \(\mathbb{P}(Z_0) = 0\) from Proposition~\ref{geo_prop}.

To evaluate the first term in (\ref{v_n_est22}), suppose that \(V_n^c \cap G_n \cap Z_0^c\) occurs. We recall that since~\(Z_0^c\) occurs, the truncated minimum passage time \(\hat{T}^{(n)}_n\) is also the passage time of the finite path~\(\hat{\pi}^{(n)}_n.\) Since \(V_n^c\) also occurs, we have that \(T_n \neq \hat{T}^{(n)}_n\) and therefore we have from Lemma~\ref{cruc} that some edge \(e \in \hat{\pi}^{(n)}_n\) has untruncated passage time \(t(e) \geq n^{\alpha}.\) But since the event \(G_n\) also occurs, every edge in the path \(\hat{\pi}^{(n)}_n\) is completely contained in the box~\(B_{n_0} = [n_0,n_0]^d,\) where \(n_0 = 8\mu \beta_1^{-1} n^{1+\epsilon_0}.\)

From the discussion in the previous paragraph we therefore have that if the event \(V_n^c \cap G_n \cap Z_0^c\) occurs, then some edge \(e \in B_{8\mu \beta_1^{-1} n^{1+\epsilon_0}}\) belongs to the geodesic~\(\hat{\pi}^{(n)}_n\) and has (untruncated) passage time \(t(e) \geq n^{\alpha}.\) Therefore
\begin{eqnarray}
\mathbb{P}\left(V_n^c \cap G_n \cap Z_0^c\right) &\leq& \mathbb{P}\left(\bigcup_{e \in B_{8\mu \beta_1^{-1} n^{1+\epsilon_0}}} \{e \in \hat{\pi}^{(n)}_n\} \cap \{t(e) \geq n^{\alpha}\} \cap Z_0^c\right) \nonumber\\
&\leq& \sum_{e \in B_{8\mu \beta_1^{-1} n^{1+\epsilon_0}}} \mathbb{P}\left(\{e \in \hat{\pi}^{(n)}_n\} \cap \{t(e) \geq n^{\alpha}\}\cap Z_0^c\right).
\;\;\;\;\;\;\;\label{v_n_est2}
\end{eqnarray}

We have the following estimate.
\begin{Lemma}\label{te_est} Fix integer \(n \geq 1\) and let \(\epsilon_0 > 0\) be a constant. For a fixed edge \(e \in B_{8\mu \beta_1^{-1} n^{1+\epsilon_0}}\) we have that
\begin{eqnarray}
\mathbb{P}\left(\{e \in \hat{\pi}^{(n)}_n\} \cap \{t(e) \geq n^{\alpha}\} \cap Z_0^c\right) \leq \frac{C}{n^{4d\alpha}},\label{t_e_box_est2}
\end{eqnarray}
for some constant \(C > 0\) independent of the choice of \(e.\) Here \(d \geq 2\) is the dimension and \(\alpha > 0\) is as in (\ref{trunc_pas}).
\end{Lemma}
\emph{Proof of Lemma~\ref{te_est}}: Fix \(e \in B_{8\mu \beta_1^{-1} n^{1+\epsilon_0}}\) and suppose that the event \[U_e := \{e \in \hat{\pi}^{(n)}_n\} \cap \{t(e) \geq n^{\alpha}\}\cap Z_0^c\] occurs. If \(x_e\) and \(y_e\) are the endvertices of \(e,\) then there are \(2d\) edge disjoint paths~\(\{P_i\}_{1 \leq i \leq 2d}\) with endvertices \(x_e\) and~\(y_e\) such that \(2d-2\) of the paths have three edges, one path is the edge \(e\) and the final path has nine edges. For example in \(d =2,\) it is easy to visualize the paths between the origin and the point \((1,0)\) as \(P_1 = ((0,0), (1,0)), P_2 = ((0,0),(0,1),(1,1),(1,0)), P_3 = ((0,0),(0,-1),(1,-1),(1,0))\) and \[P_4 = ((0,0),(-1,0),(-1,1),(-1,2),(0,2),(1,2),(2,2),(2,1),(2,0),(1,0)).\] %We remark that Cox and Durrett (1981), Kesten (1986) have also used (long) disjoint paths from the origin directly to \((n,0,\ldots,0)\) to obtain estimates on~\(T_n\) and related quantities. % with larger decay rates.

%done..WRT AS SEP PRPTY...ETC PRKMM + ETC...

For a fixed \(1 \leq i \leq 2d,\) let \(T(P_i) = \sum_{h \in P_i} t(h)\) denote the sum of passage times of edges in the path \(P_i.\) The following two properties obtain the Lemma.
\begin{equation}\label{tpi_est}
\text{If the event \(U_e\) occurs, then \(T(P_i) \geq \frac{n^{\alpha}}{2}\) for every \(1 \leq i \leq 2d.\)}
\end{equation}
For any fixed \(1 \leq i \leq 2d\) we have
\begin{equation}\label{tpi_est2}
\mathbb{P}\left(T(P_i) \geq \frac{n^{\alpha}}{2} \right)  \leq \frac{C_2}{n^{2\alpha}}
\end{equation}
for some constant \(C_2 > 0\) not depending on the choice of \(i.\)
Using properties~(\ref{tpi_est}) and~(\ref{tpi_est2}) we obtain the Lemma since
\begin{eqnarray}
\mathbb{P}\left(U_e\right) \leq \mathbb{P}\left(\bigcap_{i=1}^{2d} \left\{T(P_i) \geq \frac{n^{\alpha}}{2}\right\}\right) = \prod_{i=1}^{2d}\mathbb{P}\left(T(P_i) \geq \frac{n^{\alpha}}{2}\right) \leq \left(\frac{C_2}{n^{2\alpha}}\right)^{2d} \nonumber
\end{eqnarray}
where the equality in the middle follows since the paths \(\{P_i\}\) are edge disjoint and therefore the corresponding events are independent.

\emph{Proof of (\ref{tpi_est})}: For any finite set of edges \(A\) in \(\mathbb{Z}^d,\) let \(\hat{T}^{(n)}(A) = \sum_{e \in A} t^{(n)}(e).\) We have that
\begin{equation}\label{a_cup_b}
\hat{T}^{(n)}(A \cup B) \leq \hat{T}^{(n)}(A) + \hat{T}^{(n)}(B)
\end{equation}
with equality occurring if \(A\) and \(B\) are edge disjoint.

Let~\(Q_1\) be the subpath of the geodesic~\(\hat{\pi}^{(n)}_n\) from origin to the endvertex~\(x_e\) of the edge \(e\) and let~\(Q_2\) be the subpath of~\(\hat{\pi}^{(n)}_n\) from the endvertex~\(y_e\) of edge~\(e\) to the point~\((n,0,\ldots,0).\) The union \(Q_1 \cup Q_2 \cup\{e\} = \hat{\pi}^{(n)}_n\) and since the sets in the union are edge disjoint, we have using (\ref{a_cup_b}) that
\begin{equation} \label{t_eq_sd}
\hat{T}^{(n)}(\hat{\pi}^{(n)}_n) =  \hat{T}^{(n)}(Q_1) + {t}^{(n)}(e) + \hat{T}^{(n)}(Q_2)
\end{equation}
Since the event \(Z_0^c\) occurs, we have from Proposition~\ref{geo_prop} that the minimum truncated passage time \(\hat{T}^{(n)}_n = \hat{T}^{(n)}(\hat{\pi}^{(n)}_n).\) From (\ref{t_eq_sd}) we therefore have
\begin{eqnarray}
\hat{T}^{(n)}_n = \hat{T}^{(n)}(Q_1) + {t}^{(n)}(e) + \hat{T}^{(n)}(Q_2)= \hat{T}^{(n)}(Q_1) + n^{\alpha} + \hat{T}^{(n)}(Q_2). \label{q_split}
\end{eqnarray}
The final equality holds since \(t(e) \geq n^{\alpha}\) and so \(t^{(n)}(e) = \min(t(e),n^{\alpha}) = n^{\alpha}.\)

Suppose now that \(T(P_{i_0}) < \frac{n^{\alpha}}{2}\) for some \(1 \leq i_0 \leq 2d.\) The union of the paths \(Q_1, Q_2\) and \(P_{i_0}\) then contains a path \(Q_0\) from the origin to \((n,0,\ldots,0)\) and again using (\ref{a_cup_b}) we have
\begin{eqnarray}
\hat{T}^{(n)}(Q_0) \leq \hat{T}^{(n)}(Q_1) + \hat{T}^{(n)}(P_{i_0}) + \hat{T}^{(n)}(Q_2). \nonumber
\end{eqnarray}
Since \(t^{(n)}(h) \leq t(h)\) for all edges~\(h\) (see (\ref{trunc_pas})), we have that
\[\hat{T}^{(n)}(P_{i_0})= \sum_{h \in P_{i_0}} t^{(n)}(h) \leq \sum_{h \in P_{i_0}} t(h) = T(P_{i_0}) \leq \frac{n^{\alpha}}{2}.\]
Therefore
\begin{eqnarray}
\hat{T}^{(n)}(Q_0) \leq \hat{T}^{(n)}(Q_1) +  \frac{n^{\alpha}}{2} + \hat{T}^{(n)}(Q_2) = \hat{T}^{(n)}_n - \frac{n^{\alpha}}{2} \nonumber
\end{eqnarray}
where the final equality follows from (\ref{q_split}). But this is a contradiction since by the definition of minimum truncated passage time \(\hat{T}^{(n)}_n\) in~(\ref{t_0n_def_p}), the truncated passage time of every path with endvertices origin and \((n,0,\ldots,0)\) is at least~\(\hat{T}^{(n)}_n.\)~\(\qed\)

%Since \(e\) belongs to the geodesic \(\hat{\pi}^{(n)}_n\) and \(t(e) \geq n^{\alpha},\)

\emph{Proof of (\ref{tpi_est2})}: Fix \(1 \leq i \leq 2d\) and suppose \(T(P_i) =\sum_{h \in P_i} t(h) \geq \frac{n^{\alpha}}{2}.\) Since each path \(P_i\) has at most nine edges, we have that at least one of the edges \(h_i\) of \(P_i\) has passage time \(t(h_i) \geq \frac{n^{\alpha}}{18}.\) Therefore we have
\begin{equation}
\mathbb{P}\left(T(P_i) \geq \frac{n^{\alpha}}{2} \right) \leq \mathbb{P}\left(\bigcup_{h \in P_i} \left\{t(h) \geq \frac{n^{\alpha}}{18}\right\}\right) \leq \sum_{h \in P_i}\mathbb{P}\left(t(h) \geq \frac{n^{\alpha}}{18}\right). \label{p_i_est}
\end{equation}
Using Markov inequality we have
\[\mathbb{P}\left(t(h) \geq \frac{n^{\alpha}}{18}\right) \leq \frac{18^2}{n^{2\alpha}}\mathbb{E}t^{2}(h) \leq \frac{C_3}{n^{2\alpha}}\] for some constant~\(C_3 > 0\) not depending on \(h.\) The final estimate follows from the moment condition \((ii)\) in Section~\ref{intro}. Since there are at most nine edges in any \(P_i,\) we have from (\ref{p_i_est}) that
\[\mathbb{P}\left(T(P_i) \geq \frac{n^{\alpha}}{2} \right)  \leq \frac{9C_3}{n^{2\alpha}} = \frac{C_4}{n^{2\alpha}}\] for some constant \(C_4 > 0\) not depending on \(i.\) \(\qed\)

%From properties (\ref{tpi_est}) and (\ref{tpi_est2}) we obtain the Lemma since
%\begin{eqnarray}
%\mathbb{P}\left(U_e\right) \leq \mathbb{P}\left(\bigcap_{i=1}^{2d} \left\{T(P_i) \geq \frac{n^{\alpha}}{2}\right\}\right) = \prod_{i=1}^{2d}\mathbb{P}\left(T(P_i) \geq \frac{n^{\alpha}}{2}\right) \leq \left(\frac{C_2}{n^{2\alpha}}\right)^{2d} \nonumber
%\end{eqnarray}
%where the equality in the middle follows since the paths \(\{P_i\}\) are disjoint and therefore the corresponding events are independent. \(\qed\)

%Using properties (\ref{tpi_est}) and (\ref{tpi_est2}) we obtain the Lemma since
%\begin{eqnarray}
%\mathbb{P}\left(U_e\right) \leq \mathbb{P}\left(\bigcap_{i=1}^{2d} \left\{T(P_i) \geq \frac{n^{\alpha}}{2}\right\}\right) = \prod_{i=1}^{2d}\mathbb{P}\left(T(P_i) \geq \frac{n^{\alpha}}{2}\right) \leq \left(\frac{C_2}{n^{2\alpha}}\right)^{2d} \nonumber
%\end{eqnarray}
%where the equality in the middle follows since the paths \(\{P_i\}\) are disjoint and therefore the corresponding events are independent.

Using the estimate in the proof of Lemma~\ref{te_est}, we have the proof of Proposition~\ref{lemma_v_n}.\\
\emph{Proof of Proposition~\ref{lemma_v_n}}: We prove (\ref{v_n_est_tot}) first. Substituting the estimate (\ref{t_e_box_est2}) into (\ref{v_n_est2}) gives
\begin{eqnarray}
\mathbb{P}\left(V_n^c \cap G_n \cap Z_0^c\right) \leq \sum_{e \in B_{8\mu \beta_1^{-1} n^{1+\epsilon_0}}} \frac{C}{n^{4d\alpha}}  \leq \frac{C_1 n^{d+d\epsilon_0}}{n^{4d\alpha}} \label{v_n_est3}
\end{eqnarray}
for some constant \(C_1 > 0.\)

Setting \begin{equation}\label{eps_def}
\epsilon_0 = \frac{1}{4d} \text{ and } \alpha = \frac{1}{2} - \frac{1}{16d},
\end{equation}
we have \[4d\alpha -d - d\epsilon_0 = 2d-\frac{1}{4} -d - \frac{1}{4} = d -0.5.\]
For \(d \geq 2,\) we therefore have
\[\mathbb{P}(V_n^c \cap G_n \cap Z_0^c) \leq \frac{C_1}{n^{d-0.5}} \leq \frac{C_1}{n^{3/2}}.\] Substituting the above in (\ref{v_n_est22}) gives \[\mathbb{P}(V_n) \leq \frac{C_1}{n^{3/2}} + \frac{C_2}{n^{1+2\epsilon_0}} \leq \frac{C_3}{n^{1+2\epsilon_0}}\] for some constant \(C_3 > 0\) and for all \(n \geq 1.\) The final estimate above is true since \(1 + 2\epsilon_0 = 1 + \frac{1}{2d}  < 1 + \frac{1}{4} < \frac{3}{2}\) for all \(d \geq 2.\) This proves~(\ref{v_n_est_tot}) with \(2\epsilon_0 = \frac{1}{2d}\) as the term \(\delta\) defined in the statement of the Proposition.

To prove (\ref{w_n_est_tot}), we use the estimate (\ref{tk1k2}) of Lemma~\ref{sec_mom} to see that if \(W_n^c\) occurs then \(V^c_{n^2}\) necessarily occurs. Indeed if \(V_{n^2}\) occurs i.e., if \(T_{n^2} = \hat{T}^{(n^2)}_{n^2},\) then using (\ref{tk1k2}) we have
\[T_{n^2} = \hat{T}^{(n^2)}_{n^2} \leq \hat{T}^{(n+1)^2}_{n^2} \leq T_{n^2}.\] In other words, we have
\(\hat{T}^{(n^2)}_{n^2} = \hat{T}^{(n+1)^2}_{n^2} \) and so \(W_n\) occurs. Therefore
\begin{equation}
\mathbb{P}(W_n^c) \leq \mathbb{P}(V^c_{n^2}) \leq \frac{C_1}{n^{2+2\delta}}
\end{equation}
for some constant \(C_1 > 0,\) where the final estimate follows from~(\ref{v_n_est_tot}).

It remains to proves (\ref{t_n_eq_hat}). We have from (\ref{v_n_est_tot}) and Borel-Cantelli lemma that
\begin{equation}\label{eq_that}
\mathbb{P}\left(\liminf_n V_n \right) = 1.
\end{equation}
From (\ref{eq_that}), we have that a.e. \(\omega,\) there exists \(N_1(\omega) < \infty,\) so that for all \(n \geq N_1(\omega),\) we have \(T_n = \hat{T}^{(n)}_n.\) Therefore \[J_n := \frac{1}{n}\left(T_n - \hat{T}^{(n)}_n\right) \rightarrow 0 \text{ a.s.}\] as \(n \rightarrow \infty.\)

To see that \(\mathbb{E}J_n \rightarrow 0\) as \(n \rightarrow \infty,\) we show that \[\sup_n \mathbb{E}J_n^2 < \infty\] and this proves that \(J_n\) is uniformly integrable. Using \((a-b)^2 \leq a^2 + b^2\) for positive \(a,b,\) we have
\begin{equation}\label{j_n_est}
\mathbb{E}J_n^2 = \frac{1}{n^2} \mathbb{E}(T_n - \hat{T}^{(n)}_n)^2 \leq \frac{1}{n^2}\left(\mathbb{E}T_n^2 + \mathbb{E}\left(\hat{T}^{(n)}_n\right)^2\right) \leq 2C
\end{equation}
for some constant \(C>0.\) The final estimate follows from~(\ref{unif_int}) of Lemma~\ref{sec_mom}.\(\qed\)

\setcounter{equation}{0}
\renewcommand\theequation{\thesection.\arabic{equation}}
\section{ Variance estimates for \(T_n\) and \(\hat{T}^{(n)}_n\)} \label{var_est_sec}
In this section, we obtain variance estimates for \(\hat{T}^{(n)}_n\) and \(T_n\) needed for \(L^2\) convergence and the subsequence argument in the proof of Theorem~\ref{thm1}. %using the subsequence argument.

%We recall the definition of \(\alpha\) in truncated passage times \(t^{(k)}(e_i)\) in (\ref{trunc_pas}).

%done...CHANGE TO FINITE BOX ESTIMATES...

We have some preliminary definitions and estimates.
\subsubsection*{Boxed passage times}
To obtain the variance estimates of \(T_n\) and \(\hat{T}^{(n)}_n,\) we use the martingale difference method of Kesten~(1993) with some modifications. The proof uses the Fubini's theorem for product spaces. For simplicity and to avoid measure theoretic technicalities of infinite product spaces, we consider a ``boxed" version of the minimum passage times~\(T_n\) and~\(\hat{T}^{(n)}_n.\)

For a fixed \(\epsilon > 0,\) we have from Lemma~\ref{box_est} that the geodesic \(\hat{\pi}^{(n)}_n\) lies inside the box~\(B_{8\mu\beta^{-1}_1 n^{1+\epsilon}}\) with probability at least \(1 - \frac{C}{n^{1+2\epsilon}}\) for some constant \(C > 0.\) Here \(\mu = \sup_i \mathbb{E}t(e_i) \in (0,\infty)\) is as in~(\ref{mu_def}) and \(\beta_1 \in (0, \mu)\) is the constant in Lemma~\ref{t_pi_est} and we recall that \(B_m = [-m,m]^d\) is the box of side length~\(2m.\)

Let \(q_1,q_2,\ldots,q_N\) be the edges of the box \(B_{8\mu\beta^{-1}_1 n^{1+\epsilon}}.\) Let \(\Omega_N = \mathbb{R}^N\) and for \(\omega \in \Omega_N,\) let
\begin{equation}\label{un_def}
U_n = U_n(\omega) = \min_{\pi \subset B_{8\mu\beta^{-1}_1 n^{1+\epsilon}}} T(\pi,\omega)
\end{equation}
and
\begin{equation}\label{un_hat_def}
\hat{U}^{(n)}_n = \hat{U}^{(n)}_n(\omega) = \min_{\pi \subset B_{8\mu\beta^{-1}_1 n^{1+\epsilon}}} \hat{T}^{(n)}(\pi,\omega)
\end{equation}
be the untruncated and truncated boxed minimum passage times, respectively. For \(\omega \in \Omega_N,\) the passage times \(T(\pi,\omega) = \sum_{e \in \pi} t(e,\omega)\) and \(\hat{T}^{(n)}(\pi,\omega) = \sum_{e \in \pi} t^{(n)}(e,\omega)\) are as in (\ref{t_pi_def}) and (\ref{t_pi_p_def}), respectively. As before we suppress the dependence on \(\omega\) unless specifically mentioned.

We define~\(U_n\) and~\(\hat{U}^{(n)}_n\) on the probability space \((\Omega_N, {\cal F}_N, \mathbb{P}_{N})\) where \({\cal F}_N = \mathbb{B}(\mathbb{R}^N)\) and \(\mathbb{P}_{N}\) is the distribution of the random variables \((t(q_1),\ldots,t(q_N))\) under the measure~\(\mathbb{P}.\) We recall from Section~\ref{intro} that \(\mathbb{P}\) is the probability measure associated with the infinite sequence \((t(q_1),t(q_2),\ldots).\) For notational convenience, however, we drop the subscript from~\(\mathbb{P}_{N}\) and simply refer to it also as~\(\mathbb{P}.\) %for notational convenience.  %for notational convenience.

Fix \(\omega \in \Omega_N.\) Let \(\gamma_n = \gamma_n(\omega) \subset B_{8\mu\beta^{-1}_1 n^{1+\epsilon}}\) be the path that attains~\(U_n\) and
let \(\hat{\gamma}^{(n)}_n = \hat{\gamma}^{(n)}_n(\omega) \subset B_{8\mu\beta^{-1}_1 n^{1+\epsilon}}\) be the path that attains \(\hat{U}^{(n)}_n;\) i.e.,
\begin{equation}\label{gamma_n_def}
U_n(\omega) = T(\gamma_n,\omega) \text{ and }\hat{U}^{(n)}_n(\omega) = \hat{T}^{(n)}(\hat{\gamma}^{(n)}_n,\omega).
\end{equation}
As in Section~\ref{esti}, we refer to \(\gamma_n\) and \(\hat{\gamma}^{(n)}_n\) as geodesics and if there is more than one choice, we pick one according to a deterministic rule (see paragraph following the proof of Proposition~\ref{geo_prop}). In what follows we obtain variance estimates for \(U_n\) and \(\hat{U}^{(n)}_n\) and use those estimates to obtain variance estimates for \(T_n\) and~\(\hat{T}^{(n)}_n,\) respectively.

The following result estimates of the length (i.e. the number of edges) of the geodesics~\(\gamma_n\) and~\(\hat{\gamma}^{(n)}_n\) and is used to estimate the variances of~\(\hat{U}^{(n)}_n\) and~\(U_n.\) Let \(\#\gamma_n\) and \(\#{\hat{\gamma}^{(n)}_n}\) denote the number of edges in the respective paths. We have the following.
\begin{Lemma}\label{prop_pi_hat}We have that
\begin{equation} \label{len_pi_n2}
\mathbb{E}\left(\#{\gamma}_n\right) \leq C n
\end{equation}
and
\begin{equation} \label{len_pi_n}
\mathbb{E}\left(\#\hat{\gamma}^{(n)}_n\right) \leq C n
\end{equation}
for all \(n \geq 1\) and for some positive constant \(C.\)
\end{Lemma}
\emph{Proof of Proposition~\ref{prop_pi_hat}}: Let \(\mu = \sup_i \mathbb{E}t(e_i) \in (0,\infty)\) be as in~(\ref{mu_def}) and let \(\beta_1 \in (0, \mu)\) be the constant in Lemma~\ref{t_pi_est}. We prove (\ref{len_pi_n}) first and the proof is analogous for (\ref{len_pi_n2}). To estimate the length of \(\hat{\gamma}^{(n)}_n,\)  we have for any \(x > 0\) that
\begin{eqnarray}
\mathbb{P}(\#\hat{\gamma}^{(n)}_n \geq x) &=&  \mathbb{P}\left(\{\#\hat{\gamma}^{(n)}_n \geq x\} \cap \{\hat{T}^{(n)}(\hat{\gamma}^{(n)}_n) < \beta_1 x\}\right) \nonumber\\
&&\;\;\;\;\;\;\;+ \mathbb{P}\left(\{\#\hat{\gamma}^{(n)}_n \geq x\} \cap \{\hat{T}^{(n)}(\hat{\gamma}^{(n)}_n) \geq \beta_1 x\}\right) \nonumber\\
&\leq&\mathbb{P}\left(\{\#\hat{\gamma}^{(n)}_n \geq x\} \cap \{\hat{T}^{(n)}(\hat{\gamma}^{(n)}_n) < \beta_1 x\}\right) \nonumber\\
&&\;\;\;\;\;\;\;+ \mathbb{P}\left(\hat{T}^{(n)}(\hat{\gamma}^{(n)}_n) \geq \beta_1 x\right) \label{pi_est_eq1}
\end{eqnarray}
where \(\hat{T}^{n}(\gamma)\) is the truncated passage time of a path \(\gamma\) as defined in (\ref{t_pi_p_def}).

To estimate the second term above, we have that \[\hat{T}^{(n)}\left(\hat{\gamma}^{(n)}_n\right) \leq \sum_{i=1}^{n-1}t^{(n)}(f_i) \leq \sum_{i=1}^{n-1}t(f_i),\] where as before \(f_i\) denotes the edge between \((i-1,0,\ldots,0)\) and \((i,0,\ldots,0).\) To prove the first inequality, we argue as follows. The middle term in the above expression is the passage time of the the straight line joining the origin and \((n,0,\ldots,0)\) and this straight line is contained in the box \(B_{8\mu \beta_1^{-1} n^{1+\epsilon}}.\) The left most term is the minimum passage time among all paths with endvertices origin and \((n,0,\ldots,0)\) contained in the box \(B_{8\mu \beta_1^{-1} n^{1+\epsilon}}.\) This proves the first inequality. The second inequality holds since \(t^{(n)}(e) \leq t(e)\) for all edges \(e\) (see (\ref{trunc_pas})).

As in the proof of Lemma~\ref{lemma_a_n_est2}, we set \(X_i = t(f_i) - \mathbb{E}t(f_i)\) and have \[\mathbb{P}\left(\hat{T}^{(n)}(\hat{\gamma}^{(n)}_n) \geq \beta_1 x\right) \leq \mathbb{P}\left(\sum_{i=1}^{n-1}t(f_i)  \geq \beta_1 x\right) \leq \mathbb{P}\left(\sum_{i=1}^{n-1}X_i \geq \beta_1 x - \mu n \right),\] where the final estimate follows from the fact that \(\mathbb{E}t(f_i) \leq \mu.\)
Now for \(\beta_1 x = 8 \mu m\) and \(m \geq n\) an integer, we have \[\beta_1  x - \mu n  = 8\mu m - \mu n \geq 7 \mu m.\] Therefore \[\mathbb{P}\left(\hat{T}^{(n)}(\hat{\gamma}^{(n)}_n) \geq \beta_1 x\right) \leq \mathbb{P}\left(\sum_{i=1}^{n-1}X_i \geq 7 \mu m\right) \] and using Markov inequality, we have
\begin{eqnarray}\label{t_i_temp2}
\mathbb{P}\left(\hat{T}^{(n)}(\hat{\gamma}^{(n)}_n) \geq \beta_1 x\right) \leq \frac{\mathbb{E}\left(\sum_{i=1}^{n}X_i\right)^{2}}{(7\mu m)^{2}} \leq \frac{C_3n}{(7\mu m)^2} = \frac{C_4n}{m^2},
\end{eqnarray}
for all \(m \geq n\) and for some positive constants \(C_3\) and \(C_4.\) The second inequality above follows from (\ref{x_i_est2}).

%where Therefore (TO SEE CAREFULLY....///)
%\begin{equation}\label{t_n_pi_estp}
%\mathbb{P}\left(\hat{T}^{(n)}(\hat{\pi}^{(n)}_n) \geq 6\mu x\right) \leq \frac{constt}{(x-n)^2}
%\end{equation} using (\ref{an_prob_p}) in Lemma~\ref{lemma_a_n_est}. This estimates the second term in (\ref{pi_est_eq1}).

We now estimate the first term. Suppose now that the event in the first term of (\ref{pi_est_eq1}) occurs with \(\beta_1 x = 8\mu m,\) for some \(m \geq n.\) This implies that there exists a path \(\pi (= \hat{\gamma}^{(n)}_n)\) containing \(r \geq x = \frac{8\mu}{\beta_1}m\) edges with truncated passage time \(\hat{T}^{(n)}(\pi) < \beta_1 x \leq \beta_1 r.\) In particular, the event \({E}_{m}(n)\) defined in~(\ref{e_k_hat_def}) occurs and so for \(m \geq n\) we have using (\ref{a_0k_p}) that
\begin{equation}\label{t_n_pi_estp2}
\mathbb{P}\left(\{\#{\hat{\gamma}^{(n)}_n} \geq 8\mu\beta_1^{-1}m\} \cap \{\hat{T}^{(n)}(\hat{\gamma}^{(n)}_n) < 8 \mu m\}\right) \leq \mathbb{P}({E}_{m}(n)) \leq C_2e^{-\beta_2 m}
\end{equation}
where \(C_2, \beta_2 > 0\) are as in (\ref{a_0k_p}).

%Since \(\beta_1 < \mu\) by choice, we have that \(8\mu\beta_1^{-1}m \geq 8\mu m.\) (to see carefully range etc...) \(\frac{constt}{(8\mu\beta_1^{-1}m-n)^2} \leq \)

Substituting (\ref{t_i_temp2}) and (\ref{t_n_pi_estp2}) into (\ref{pi_est_eq1}) gives
\begin{eqnarray}
\mathbb{P}(\#\hat{\gamma}^{(n)}_n \geq 8\mu\beta_1^{-1}m) \leq \frac{C_4 n}{m^2} + C_2e^{-\beta_2 m}.\label{pi_est_eq2}
\end{eqnarray}
for all \(m \geq n.\) We therefore have
\begin{eqnarray}
\frac{1}{8\mu\beta_1^{-1}}\mathbb{E}(\#{\hat{\gamma}^{(n)}}_n) &\leq& \sum_{m \geq 0} \mathbb{P}(\#{\hat{\pi}^{(n)}_n} \geq 8\mu \beta_1^{-1}m) \nonumber\\
&=& (\sum_{m = 0}^{n-1}  + \sum_{m  \geq n}) \mathbb{P}(\#{\hat{\gamma}^{(n)}_n} \geq 8\mu \beta_1^{-1}m)  \nonumber\\
&\leq& n + \sum_{m  \geq n} \mathbb{P}(\#{\hat{\gamma}^{(n)}_n} \geq 8\mu \beta_1^{-1}m) \nonumber\\
&\leq& n + \sum_{m \geq n} \frac{C_4 n}{m^2} + \sum_{m  \geq n} C_2e^{-\beta_2 m}. \nonumber
\end{eqnarray}
We have \[\sum_{m \geq n} \frac{C_4 n}{m^2} \leq C_5 \text{ and }\sum_{m \geq n} C_2e^{-\beta_2 m} \leq C_6 e^{-\beta_2 n}\] for some positive constants \(C_5\) and \(C_6.\) Therefore
\[\frac{1}{8\mu\beta_1^{-1}}\mathbb{E}(\#{\hat{\gamma}^{(n)}}_n) \leq n +  C_5  + C_6 e^{-\beta_2 n} \leq 2n\] for all \(n\) large.\(\qed\)

\subsection{Variance Estimates}
We recall that the boxed passage time \(\hat{U}^{(n)}_n\) is the minimum of all passage times of paths contained in the box \(B_m,\) where \(m = 8\mu\beta_1^{-1}n^{1+\epsilon}.\) Here \(\mu = \sup_i \mathbb{E}t(e_i) \in (0,\infty)\) is as in~(\ref{mu_def}) and \(\beta_1 \in (0, \mu)\) is the constant in Lemma~\ref{t_pi_est}. We have the following result regarding the variance of the boxed passage times.
\begin{Lemma}\label{eun} We have that
\begin{equation}\label{eun_est}
\mathbb{E}(\hat{U}^{(n)}_n - \mathbb{E}\hat{U}^{(n)}_n)^2 \leq C_1 n
\end{equation}
for all \(n \geq 1\) and some constant \(C_1 > 0.\) Setting \(\epsilon = 3\) and using (\ref{eun_est}), we also have that
\begin{equation}\label{etn_est}
\mathbb{E}(\hat{T}^{(n)}_n - \mathbb{E}\hat{T}^{(n)}_n)^2 \leq C_2 n
\end{equation}
for all \(n \geq 1\) and some constant \(C_2 > 0.\)
\end{Lemma}

Set \({\cal F}_0 = \{ \emptyset, \Omega\}\) and for integer \( 1 \leq i \leq N,\) set \[{\cal F}_i = \sigma({t}^{(n)}(q_l) : 1 \leq l \leq i)\] to be the sigma field generated by the truncated passage time of the edges~\(\{q_l\}_{1 \leq l \leq i}.\)
For \(1 \leq l \leq N,\) let
\begin{equation}\label{x_l_def}
X_{l} := \mathbb{E}(\hat{U}^{(n)}_n|{\cal F}_l) - \mathbb{E} (\hat{U}^{(n)}_n|{\cal F}_{l-1}).
\end{equation}
There is a finite path \(\hat{\gamma}^{(n)}_n \subset B_{8\mu\beta^{-1}_1 n^{1+\epsilon}}\) whose passage time is~\(\hat{U}^{(n)}_n.\)
The following estimate regarding~\(X^2_l\) is used in the proof of Lemma~\ref{eun}.
\begin{Lemma}\label{x_l_est} For \(1 \leq l \leq N,\) we have
\begin{equation}
\mathbb{E}(X_l^2|{\cal F}_{l-1}) \leq C_1 \mathbb{P}(q_l \in \hat{\gamma}^{(n)}_n | {\cal F}_{l-1}) \text{ a.s. }\label{eq1}
\end{equation}
for some positive constant \(C_1\) not depending on \(l\) or \(n.\)
\end{Lemma}

\emph{Proof of Lemma~\ref{x_l_est}}: Fix \(1 \leq l \leq N.\) For \(1 \leq j \leq N,\) let \({\nu}_j(.)\) denote the probability measure associated with the random vector \[({t}^{(n)}(q_j),{t}^{(n)}(q_{j+1}),\ldots,t^{(n)}(q_N)).\] Let \(\sigma = (\sigma_1,\sigma_2,\ldots, \sigma_N)\) and \(\omega = (\omega_1,\omega_2,\ldots,\omega_N) \in \Omega_N\) and define
\begin{equation}\label{omeg_a_def}
\omega_a = (\omega_1,\omega_2,\ldots,\omega_{l},\sigma_{l+1},\sigma_{l+2},\ldots,\sigma_N)
\end{equation}
and
\begin{equation}\label{omeg_b_def}
\omega_b = (\omega_1,\omega_2,\ldots,\omega_{l-1},\sigma_{l},\sigma_{l+2},\ldots,\sigma_N).
\end{equation}
If \(l = N,\) then \(\omega_a = (\omega_1,\ldots,\omega_N)\) and if \(l = 1,\) then \(\omega_b = (\sigma_1,\ldots,\sigma_N).\)
In the notation of Kesten (1993), \(\omega_a = [\omega,\sigma]_l\) and \(\omega_b = [\omega,\sigma]_{l-1}.\)

We introduce the following temporary notation. For any fixed path \(\gamma,\) we let \(a(\gamma) = \hat{T}^{(n)}(\gamma,\omega_a)\) and \(b(\gamma) = \hat{T}^{(n)}(\gamma,\omega_b)\) be the truncated passage times of the path \(\gamma\) for the configurations \(\omega_a\) and \(\omega_b,\) respectively. Let \(\hat{U}^{(n)}_n(\omega_a) = \hat{T}^{(n)}(\gamma_a,\omega_a) = a(\gamma_a)\) be the value of the minimum boxed passage time~\(\hat{U}^{(n)}_n\) for realization \(\omega_a\) as defined in (\ref{un_def}) and~\(\gamma_a\) be the geodesic with the minimum boxed passage time.

Throughout this proof we only consider paths that are completely contained in the box~\(B_m\) where \(m = 8\mu\beta_1^{-1}n^{1+\epsilon}.\)
We have using Fubini's theorem that
\begin{equation}\label{x_l_alt}
X_l = X_l(\omega) = \int \nu_l(d\sigma) W_{l}\;\;\;\text{  a.e. }\omega
\end{equation}
where
\begin{equation}\label{w_nl}
W_{l} := \hat{U}^{(n)}_n(\omega_a) - \hat{U}^{(n)}_n(\omega_b) = a(\gamma_a) - b(\gamma_b)
\end{equation}
We find a good estimate on~\(W_{l}\) as follows. If the edge \(q_l\) does not belong to the geodesic~\(\gamma_a\) and does not belong to~\(\gamma_b,\) then \(\gamma_a = \gamma_b\) and \(a(\gamma_a) = b(\gamma_b).\) Therefore,
\begin{eqnarray}
|a(\gamma_a) - b(\gamma_b)| = |a(\gamma_a) - b(\gamma_b)|\ind\left(\{q_l \in \gamma_a\} \cup \{q_l \in \gamma_b\} \right). \label{first_est_ab}
\end{eqnarray}
For any fixed path \(\gamma,\) we have that \(a(\gamma) = b(\gamma)\) if \(q_l \notin \gamma\) since the passage times of any path not containing the edge~\(q_l\) are the same in both configurations. Similarly if \(q_l \in \gamma,\) then
\begin{equation}\label{a_pi_dif}
|a(\gamma)-b(\gamma)| = |Y_a - Y_b|,
\end{equation}
where  \(Y_a = t^{(n)}\left(q_l,\omega_a\right)\) and \(Y_b =  t^{(n)}\left(q_l,\omega_b\right)\) denote the truncated passage times of the edge \(q_l\) in the respective configurations \(\omega_a\) and \(\omega_b.\)
In particular,
\begin{equation}\label{a_pi_b_pi}
|a(\gamma)-b(\gamma)| \leq |Y_a - Y_b|
\end{equation}
for any finite path \(\gamma.\) Writing \(a(\gamma) \leq b(\gamma) + |Y_a-Y_b|\) and taking minimum over all finite paths contained in the box~\(B_{8\mu\beta_1^{-1}n^{1+\epsilon}}\) and using the fact that \(\gamma_a\) and \(\gamma_b\) are the minimum boxed passage times for~\(\omega_a\) and~\(\omega_b,\) respectively, we have that \(a(\gamma_a) \leq b(\gamma_b)  + |Y_a-Y_b|.\) Similarly, using the other inequality in (\ref{a_pi_b_pi}) and taking minimum again, we get \(b(\gamma_b) \leq a(\gamma_a) + |Y_a-Y_b|.\) Thus we have \[|a(\gamma_a) - b(\gamma_b)| \leq |Y_a-Y_b|.\]

Substituting the above in (\ref{first_est_ab}) gives
\begin{equation}\label{a_pi_est2}
|a(\gamma_a) - b(\gamma_b)| \leq |Y_a-Y_b|\ind(A_a \cup A_b),
\end{equation}
where \(A_a = \{q_l \in \gamma_a\}\) and \(A_b = \{q_l \in \gamma_b\}.\)

We improve the estimate (\ref{a_pi_est2}) to obtain
\begin{equation}\label{a_pi_est3}
|a(\pi_a) - b(\pi_b)| \leq Y_b \ind(Y_b > Y_a)\ind(A_a) + Y_a \ind (Y_a \geq Y_b) \ind (A_b).
\end{equation}
\emph{Proof of (\ref{a_pi_est3})}: We assume \(Y_a \geq Y_b\) and obtain the second term. An analogous analysis holds for the first term. If \(Y_a = t^{(n)}(q_l,\omega_a) \geq t^{(n)}(q_l,\omega_b) = Y_b\) then \[|Y_a - Y_b|  = Y_a - Y_b \leq Y_a.\] Moreover if \(\ind(A_a \cup A_b) = 1\) also holds, then necessarily \(A_b = \{q_l \in \gamma_b\}\) occurs. This is true if \(Y_a = Y_b\) since then both the configurations are identical. In what follows we assume that \(Y_a > Y_b\) strictly. Roughly speaking, lowering the passage time of the edge \(q_l\) from~\(Y_a\) in the ``old" configuration~\(\omega_a\) (with geodesic \(\gamma_a\)) to \(Y_b\) in the ``new" configuration~\(\omega_b\) (with geodesic~\(\gamma_b\)), improves the chances of~\(q_l\) belonging to the geodesic~\(\gamma_b\) of the new configuration.

More formally, suppose the event \(A_a \cup A_b\) holds. If \(A_a\) does not occur, then the event \(A_b\) necessarily occurs and we are done. Suppose now that the event \(A_a\) occurs. We then have that the edge \(q_l\) belongs to the geodesic~\(\pi_a\) of the old configuration~\(\omega_a.\) Also the passage time \(Y_a\) of~\(q_l\) in \(\omega_a\) is strictly larger than the passage time~\(Y_b\) in the new configuration~\(\omega_b\) and the passage time of every other edge remains the same in both configurations. We therefore have that
\begin{eqnarray}
\hat{T}^{(n)}\left(\gamma_a,\omega_b\right) = \sum_{e \in \gamma_a} t^{(n)}\left(e,\omega_b\right) < \sum_{e \in \gamma_a} {t}^{(n)}\left(e,\omega_a\right)  = \hat{T}^{(n)}\left(\gamma_a,\omega_a\right). \label{tn_temp_ab1}
\end{eqnarray}
In words, the truncated passage time of the path \(\gamma_a\) in the new configuration~\(\omega_b\) is strictly less than the passage time in the old configuration~\(\omega_a.\)

The final term in (\ref{tn_temp_ab1}) is the minimum boxed passage time \(\hat{U}^{(n)}_n(\omega_a)\) in the old configuration \(\omega_a.\) In particular, this means that the minimum boxed passage time \(\hat{U}^{(n)}_n(\omega_b)\) in the new configuration \(\omega_b\) is strictly less than the minimum boxed passage time \(\hat{U}^{(n)}_n(\omega_a)\) in the old configuration~\(\omega_a.\) Therefore
\begin{eqnarray}
\hat{T}^{(n)}\left(\gamma_b,\omega_b\right) = \hat{U}^{(n)}_n(\omega_b) < \hat{U}^{(n)}_n\left(\omega_a\right) \leq \hat{T}^{(n)}\left(\gamma_b,\omega_a\right) \nonumber
\end{eqnarray}
where the final estimate follows from the definition of minimum boxed passage time in (\ref{un_def}).

The above estimate necessarily means that \(q_l \in \gamma_b\) since for any fixed path \(\gamma\) not containing the edge \(q_l\) we have that the passage time of~\(\gamma\) in both the configurations remains the same; i.e.,
\begin{eqnarray}
\hat{T}^{(n)}\left(\gamma,\omega_a\right) = \sum_{e \in \gamma} {t}^{(n)}\left(e,\omega_a\right) = \sum_{e \in \pi} {t}^{(n)}\left(e,\omega_b\right) = \hat{T}^{(n)}\left(\gamma,\omega_b\right). \nonumber
\end{eqnarray}
This proves the second term in~(\ref{a_pi_est3}). The proof for the first term is analogous.~\(\qed\)

Squaring both sides of (\ref{a_pi_est3}) gives
\begin{equation}\label{a_pi_est22}
|a(\gamma_a) - b(\gamma_b)|^2  \leq Y_b^2\ind(Y_b > Y_a)\ind(A_a)  + Y_a^2 \ind (Y_a \geq Y_b) \ind (A_b).
\end{equation}

Using Cauchy-Schwarz inequality in~(\ref{x_l_alt}), we first for a.e. \(\omega\) that
\begin{equation}
X_l^2 \leq \int \nu_l(d\sigma) W^2_{l} \leq \int \nu_l(d\sigma) |a(\gamma_a) - b(\gamma_b)|^2 \leq 2Z_1  + 2Z_2,\label{x_l_2_def}
\end{equation}
where \[Z_1 = \int \nu_l(d\sigma) Y_b^2 \ind(Y_b > Y_a)\ind(A_a) \text{ and } Z_2 = \int \nu_l(d\sigma) Y_a^2 \ind (Y_a \geq Y_b) \ind (A_b).\]

Taking conditional expectation w.r.t \({\cal F}_{l-1}\) we have for a.e. \(\omega\) that
\begin{equation}\label{x_l_3_def}
\mathbb{E}(X_l^2 | {\cal F}_{l-1}) \leq 2 \mathbb{E}(Z_1 | {\cal F}_{l-1}) + 2 \mathbb{E}(Z_2 | {\cal F}_{l-1}).
\end{equation}
We estimate the two terms above separately. By Fubini's theorem we first have
for a.e. \(\omega\) that
\begin{equation}\label{cond_exp_al}
\int \nu_l(d\sigma) \ind(A_a) = S_a \text{ and } \int \nu_l(d\sigma) \ind(A_b) = S_b.
\end{equation}
where \[S_a = S_a(\omega_1,\ldots,\omega_l) = \mathbb{P}\left(q_l \in \hat{\gamma}^{(n)}_n | {\cal F}_l\right) \] and \[S_b = S_b(\omega_1,\ldots, \omega_{l-1}) = \mathbb{P}\left(q_l \in \hat{\gamma}^{(n)}_n | {\cal F}_{l-1}\right).\] By property of conditional expectation, we also have that
\begin{equation}\label{f_1_est}
\mathbb{E}(S_a | {\cal F}_{l-1}) = S_b \text{ a.e. } \omega.
\end{equation}

To evaluate \(Z_1\) we recall that the configurations \(\omega_a\) and \(\omega_b\) defined in (\ref{omeg_a_def}) and (\ref{omeg_b_def}), respectively. The edge passage time \(Y_b = t^{(n)}(q_l,\omega_b)\) depends only on the variable \(\sigma_l\) and the term \(\ind(A_a) = \ind(A_a)(\omega_a)\) does not depend on \(\sigma_l\) since \(\omega_a = (\omega_1,\ldots,\omega_l,\sigma_{l+1},\ldots,\sigma_N).\) Letting \(\mu_l\) denote the probability measure associated with \(t^{(n)}(q_l,\sigma_l),\) we therefore have for a.e. \(\omega\)
that
\begin{eqnarray}
Z_1 &=& \int \nu_l(d\sigma) Y_b^2 \ind(Y_b > Y_a)\ind(A_a) \nonumber\\
&\leq& \int \nu_l(d\sigma)  Y_b^2(\sigma_l) \ind(A_a)\nonumber\\
&=& \int \mu_l(d\sigma_l)  Y_b^2(\sigma_l) \int \nu_{l+1} (d\sigma)\ind(A_a). \label{z1_temp_est}
\end{eqnarray}
Using (\ref{cond_exp_al}), we have that the second term in (\ref{z1_temp_est}) is
\[\int \nu_{l+1} (d\sigma)\ind(A_a) = \int \nu_{l} (d\sigma)\ind(A_a) = S_a(\omega_1,\ldots,\omega_l).\] The middle equality is true since \(\ind(A_a)(\omega_a)\) does not depend on the variable~\(\sigma_l.\)

The first term in (\ref{z1_temp_est}) is
\[\int \mu_l(d\sigma_l)  Y_b^2(\sigma_l) = \mathbb{E}\left(t^{(n)}(q_l)\right)^{2}  \leq \mathbb{E}t^2(q_l) \leq \sup_i \mathbb{E}t^2(q_i) \leq C \]
for some constant \(C > 0.\) The first inequality follows from the fact that \(t^{(n)}(q_l) \leq t(q_l)\) for any \(n \geq 1\) (see (\ref{trunc_pas})) and the final inequality is true by the moment condition \((ii)\) in Section~\ref{intro}.
Thus
\begin{equation}\label{z_1_est1}
Z_1 \leq C S_a(\omega_1,\ldots,\omega_l)
\end{equation}
and taking conditional expectations and using (\ref{f_1_est}) we have a.e. \(\omega\) that
\begin{equation}\label{z_1_est2}
\mathbb{E}(Z_1 | {\cal F}_{l-1}) \leq C S_b = C \mathbb{P}(q_l \in \hat{\pi}^{(n)}_n | {\cal F}_{l-1}).
\end{equation}

%//....CHCK BELOW ONWARDS...22 10

%//....CHCK ALSO END CONDITIONS I.E. L = 1 AND L = N ETC...PRKVMM+ETC...

An analogous argument holds for \(Z_2\) in (\ref{x_l_2_def}); we note that \(Y_a = Y_a(\omega_l) = t^{(n)}(q_l,\omega_l)\) is independent of \(\ind(A_b) = \ind(A_b)(\omega_b)\) since \(\omega_b = (\omega_1,\ldots,\omega_{l},\sigma_{l+1},\ldots,\sigma_N).\)
Therefore a.e. \(\omega\) we have
\begin{eqnarray}
Z_2 &=& \int \nu_l(d\sigma) Y_a^2 \ind(Y_a \geq Y_b)\ind(A_b) \nonumber\\
&\leq& \int \nu_l(d\sigma)  Y_a^2(\omega_l) \ind(A_b)\nonumber\\
&=& Y_a^2(\omega_l) \int \nu_{l} (d\sigma)\ind(A_b)\nonumber\\
&=& Y_a^2(\omega_l) S_b(\omega_1,\ldots,\omega_{l-1}) \nonumber
\end{eqnarray}
where the last equality follows from (\ref{cond_exp_al}).

Again taking conditional expectation w.r.t \({\cal F}_{l-1}\) gives a.e. \(\omega\) that
\begin{eqnarray}
\mathbb{E}(Z_2 | {\cal F}_{l-1}) &\leq& \int \nu_l(d\omega) Y_a^2(\omega_l) S_b(\omega_1,\ldots,\omega_{l-1})  \nonumber\\
&=& S_b(\omega_1,\ldots,\omega_{l-1}) \int \nu_l(d\omega) Y_a^2(\omega_l).\label{z_2_est1}
\end{eqnarray}
The first equality follows since \(S_b\) depends only on \(\{\omega_i\}_{1 \leq i \leq l-1}\)
and the term under the integral is estimated as
\[\int \nu_l(d\omega) Y_a^2(\omega_l)  = \mathbb{E}\left(t^{(n)}(q_l)\right)^2  \leq \mathbb{E}t^2(q_l)  \leq C_3\]
for some constant \(C_3 > 0.\) The first inequality above is true since \(t^{(n)}(q_l) \leq t(q_l)\) for any fixed \(n \geq 1.\)  The final estimate follows from the moment condition~\((ii)\) in Section~\ref{intro}.

Using (\ref{z_1_est2}) and (\ref{z_2_est1}) in (\ref{x_l_3_def}) gives the required estimate in (\ref{eq1}). \(\qed\)

%Using (\ref{pi_a_def}) and (\ref{pi_b_def}) again in the above gives \[|\hat{T}^{(n)}_n([\omega,\sigma]_l) - \hat{T}^{(n)}_n([\omega,\sigma]_{l-1})| \leq n^{\alpha} \left(\ind(e_l \in \hat{\pi}^{(n)}_n([\omega,\sigma]_l)) + \ind(e_l \in \hat{\pi}{^(n)}_n([\omega,\sigma]_{l-1}))\right).\] By Cauchy-Schwarz inequality, we have the following inequalities hold a.s.;
%\begin{eqnarray}
%X_l^2 &\leq& \int {\nu}_l(d\sigma)|\hat{T}^{(n)}_n([\omega,\sigma]_l) - \hat{T}^{(n)}_n([\omega,\sigma]_{l-1})|^2 \nonumber\\
%&\leq& 2n^{2\alpha} \int {\nu}_l(d\sigma)\left(\ind(e_l \in {\pi}_n([\omega,\sigma]_l)) + \ind(e_l \in {\pi}_n([\omega,\sigma]_{l-1}))\right) \nonumber\\
%&=& 2n^{2\alpha} \left(\mathbb{P}(e_l \in {\pi}_n | {\cal F}_l) + \mathbb{P}(e_l \in {\pi}_n | {\cal F}_{l-1})\right). \nonumber
%\end{eqnarray}
%This proves (\ref{eq1}).

%done...WRITE ABOUT THE BOXED ETC... IN PROOF BELOW +ETC AND MOR PROVMM +ETC..\\
\emph{Proof of Lemma~\ref{eun}}: We first prove (\ref{eun_est}). We recall that \(N\) denotes the number of edges in the box~\(B_{8\mu \beta_1^{-1}n^{1+\epsilon}}.\)
From the definition of \(X_l, 1 \leq l \leq N\) in (\ref{x_l_def}) we have that
\begin{equation}\label{tel_scop}
\sum_{l=1}^{N} X_l = \mathbb{E} (\hat{U}^{(n)}_n | {\cal F}_N) - \mathbb{E}(\hat{U}^{(n)}_n | {\cal F}_0) = \hat{U}^{(n)}_n - \mathbb{E}\hat{U}^{(n)}_n
\end{equation}
since \(\hat{U}^{(n)}_n\) defined in (\ref{un_hat_def}) is \({\cal F}_N-\)measurable.

By the martingale property, we have that
\begin{equation}\label{y_m_est}
\mathbb{E}\left(\sum_{l=1}^{N}X_l\right)^{2} = \sum_{l=1}^{N} \mathbb{E}X_l^2.
\end{equation}
\emph{Proof of (\ref{y_m_est})}: We have that
\begin{equation}\label{y_m_est2}
\mathbb{E}\left(\sum_{l=1}^{N}X_l\right)^2 = \sum_{l=1}^{N}\mathbb{E}X_l^2 + 2\sum_{i=1}^{N-1}\sum_{j=i+1}^{N}\mathbb{E}X_iX_j.
\end{equation}
If \(j \geq i+1,\) we have \[\mathbb{E}(X_iX_j) = \mathbb{E}\mathbb{E}(X_iX_j | {\cal F}_i) = \mathbb{E}\left(X_i\mathbb{E}(X_j|{\cal F}_i)\right).\] The final equality follows from the definition of \(X_l\) in~(\ref{x_l_def}) where we see that \(X_l\) is \({\cal F}_l-\) measurable since \({\cal F}_{l-1} \subseteq {\cal F}_l\) for all \(1 \leq l \leq N.\) We therefore have that if \(j \geq i+1,\) then
\begin{eqnarray}
\mathbb{E}(X_j | {\cal F}_i) &=& \mathbb{E}\left(\mathbb{E}(\hat{U}^{(n)}_n | {\cal F}_j)|{\cal F}_i\right) - \mathbb{E}\left(\mathbb{E}(\hat{U}^{(n)}_n | {\cal F}_{j-1})|{\cal F}_i\right)\nonumber\\
&=& \mathbb{E}(\hat{U}^{(n)}_n | {\cal F}_i) - \mathbb{E}(\hat{U}^{(n)}_n | {\cal F}_i) \nonumber\\
&=& 0. \nonumber
\end{eqnarray}
Substituting the above into (\ref{y_m_est2}) gives (\ref{y_m_est}).~\(\qed\)

From (\ref{y_m_est}) and (\ref{tel_scop}), we thus have
\begin{equation}\nonumber
\mathbb{E}(\hat{U}^{(n)}_n -\mathbb{E}\hat{U}^{(n)}_{n})^2  = \sum_{l=1}^{N}\mathbb{E}X_l^2.
\end{equation}
Using the estimate (\ref{eq1}) for \(X_l^2,\) we have for \(1 \leq l \leq N\) that
\begin{eqnarray}
\mathbb{E}X_l^2 = \mathbb{E} \left(\mathbb{E}(X_l^2 | {\cal F}_{l-1}) \right) \leq C_1 \mathbb{E} \left(\mathbb{P}(q_l \in \hat{\gamma}^{(n)}_n |{\cal }F_{l-1})\right) = C_1 \mathbb{P}(q_l \in \hat{\gamma}^{(n)}_n) \nonumber
\end{eqnarray}
for some constant \(C_1 >0.\)
Thus
\begin{equation}\label{fin_est_t_n}
\mathbb{E}(\hat{U}^{(n)}_n -\mathbb{E}\hat{U}^{(n)}_{n})^2 \leq C_1 \sum_{l=1}^{N} \mathbb{P}(q_l \in \hat{\gamma}^{(n)}_n) = C_1 \mathbb{E}\sum_{l=1}^{N} \ind(q_l \in \hat{\gamma}^{(n)}_n)
= C_1\mathbb{E}(\#{\hat{\gamma}^{(n)}_n})
\end{equation}
where \(\ind(.)\) refers to the indicator function and \(\#{\hat{\gamma}^{(n)}_n}\) refers to the number of edges in the path~\({\hat{\gamma}^{(n)}_n}.\) %and the first inequality follows from~(\ref{eq1}).
In Lemma~\ref{prop_pi_hat} we have estimated the number of edges in \(\hat{\gamma}^{(n)}_n\) to be at most a constant multiple of \(n\) (see estimate (\ref{len_pi_n}) of Lemma~\ref{prop_pi_hat}). Substituting~(\ref{len_pi_n}) into the above expression, we have the estimate (\ref{eun_est}) of the Lemma.

To obtain estimate (\ref{etn_est}), we first have a small observation.
For any two random variables \(X\) and \(Y\) we have \[var (X+Y) = \mathbb{E}\left( (X - \mathbb{E}X) + (Y - \mathbb{E}Y)\right)^2.\] Using \((a+b)^2 \leq 2(a^2 + b^2)\) for any two real numbers \(a\) and \(b,\) we have
\(var(X+Y) \leq 2(var(X) + var(Y)).\) Letting \(X = \hat{T}^{(n)}_n\) and \(Y = \hat{U}^{(n)}_n,\) we have
\begin{equation}\label{tn_un_est}
var\left(\hat{T}^{(n)}_n \right) \leq 2 var\left(\hat{T}^{(n)}_n  - \hat{U}^{(n)}_n\right)  + 2 var\left(\hat{U}^{(n)}_n \right).
\end{equation}
The second term in (\ref{tn_un_est}) is at most \(2C_1 n\) using (\ref{eun_est}). To estimate the first term we recall that the term \(\hat{U}^{(n)}_n\) is the minimum passage time of all paths contained in the box \(B_{8\mu \beta_1^{-1} n^{1+\epsilon}}.\) Setting \(\epsilon = 3\) we have
\begin{equation}\label{un_est}
var\left(\hat{T}^{(n)}_n - \hat{U}^{(n)}_n \right) \leq \frac{C_2}{n^4}
\end{equation}
for some constant \(C_2 > 0\) and all \(n \geq 1.\) From (\ref{tn_un_est}), we therefore have that \[var\left(\hat{T}^{(n)}_n \right) \leq 2 \frac{C_2}{n^4}  + 2C_1n \leq C_3n\] for some constant \(C_3 >0\) and for all \(n \geq 1.\) This proves (\ref{etn_est}) of the Lemma.

\emph{Proof of (\ref{un_est})}: Let
\begin{equation}\nonumber
\hat{H}_n = \left\{\hat{\pi}^{(n)}_n \subseteq B_{8\mu \beta_1^{-1} n^{1+\epsilon}}\right\}
\end{equation}
be the event that every edge of the geodesic~\(\hat{\pi}^{(n)}_n \) for the truncated passage time \(\hat{T}^{(n)}_n\) is contained in the box~\(B_{8\mu \beta_1^{-1} n^{1+\epsilon}}.\)
From the definition of the passage times \(\hat{U}^{(n)}_n\) and \(\hat{T}^{(n)}_n\) we have that
\begin{equation}\label{tu_eq}
\hat{U}^{(n)}_n \ind(\hat{H}_n) = \hat{T}^{(n)}_n \ind(\hat{H}_n).
\end{equation}

Letting \(\hat{Z}_n = \hat{T}^{(n)}_n  - \hat{U}^{(n)}_n,\) we have
\begin{equation}\label{sec_term_est}
var(\hat{Z}_n) \leq \mathbb{E}(\hat{Z}_n)^2 = \mathbb{E}(\hat{Z}_n)^2\ind(\hat{H}_n^c).
\end{equation}
We have that \[|\hat{Z}_n| \leq \hat{T}^{(n)}_n + \hat{U}^{(n)}_n \leq 2\sum_{i=1}^{n}t^{(n)}(f_i) \leq 2n^{1+\alpha}\] where \(\alpha > 0\) is as in (\ref{trunc_pas}). The second inequality is true because the terms~\(\hat{T}^{(n)}_n\) and \(\hat{U}^{(n)}_n\) are each no more than the truncated passage time of the straight line with endvertices origin and~\((n,0,\ldots,0).\) The final estimate is true since \(t^{(n)}(e) \leq n^{\alpha}\) for any edge \(e,\) by definition (see (\ref{trunc_pas})).
Thus from~(\ref{sec_term_est}) we have
\[var(\hat{Z}_n) \leq 4n^{2+2\alpha}\mathbb{P}(\hat{H}_n^c) \leq n^{2+2\alpha}\frac{C_2}{n^{1+2\epsilon}} = \frac{C_2}{n^{2\epsilon-2\alpha-1}}\]  for some constant \(C_2 > 0\) where the second inequality follows from (\ref{box_pi}) of Lemma~\ref{box_est}. From the relation (\ref{eps_def}), we have that \(\alpha < \frac{1}{2}.\) Setting \(\epsilon = 3,\) we therefore have \(2\epsilon - 2\alpha -1 > 2\epsilon - 2 = 4.\) This proves (\ref{un_est}).~\(\qed\)

%if
%the minimum truncated passage time \(\hat{T}^{(n)}_n\) is no more than  the passage time \(\sum_{i=1}^{n} t^{(n)}(f_i)\) of the straight line with endvertices origin and~\((n,0,\ldots,0).\)

%From (\ref{box_est}), we have that \(\mathbb{P}(G^c_n) \leq \frac{1}{n^{1+2\epsilon}}\) and so from Borel-Cantelli Lemma, we have that \(\mathbb{P}(G^c_n \text{ i. o.}) = 0.\) Therefore \(\ind(G_n^c) \rightarrow 0\) a.s.
%Also from estimate () of Lemma~\ref{sec_mom} we have \(Z_n \leq \sum_{i=1}^{n}t(f_i)\) and

We have an analogous result for the variance of the untruncated random variable \(U_n\) and \(T_n\) with some minor differences. %we have an analogous result.
\begin{Lemma}\label{eun2} We have that
\begin{equation}\label{eun_est2}
\mathbb{E}(U_n - \mathbb{E}U_n)^2 \leq C_1 n
\end{equation}
for all \(n \geq 1\) and some constant \(C_1 > 0.\) If the uniform integrability condition~\((ii)(a)\) of Section~\ref{intro} holds, then setting \(\epsilon = 3\) and using (\ref{eun_est2}) we have that
\begin{equation}\label{etn_est2}
\mathbb{E}\left(\frac{T_n}{n} - \frac{\mathbb{E}T_n}{n}\right)^2 \longrightarrow 0
\end{equation}
as \(n \rightarrow \infty.\) If the condition~\((ii)(b)\) of Section~\ref{intro} holds, then setting \(\epsilon = \frac{3}{p-2}+1\) and using (\ref{eun_est2}) we have
\begin{equation}\label{etn_est3}
\mathbb{E}(T_n-\mathbb{E}T_n)^2 \leq Cn
\end{equation}
for some constant \(C > 0\) and for all \(n \geq 1.\)
\end{Lemma}
Thus if the passage times have a bounded \(p^{th}\) moment for some \(p > 2,\) we have that the variance of the passage time \(T_n\) is at most a constant multiple of \(n.\) %The estimate (\ref{etn_est2}) is slightly weaker than the corresponding estimate~(\ref{etn_est}) for the truncated minimum passage time.

\emph{Proof of Lemma~\ref{eun2}}: The proof of (\ref{eun_est2}) is analogous to (\ref{eun_est}). We use the notation in the proof of Lemma~\ref{eun}. We set \({\cal G}_0 = \{ \emptyset, \Omega\}\) and for each \( 1 \leq i \leq N,\) set \[{\cal G}_i = \sigma({t}(q_l) : 1 \leq l \leq i)\] to be the sigma field generated by the untruncated passage time of the edges \(\{q_l\}_{1 \leq l \leq i}.\) For \(1 \leq l \leq N,\) we reuse notation of Lemma~\ref{x_l_est} and let
\begin{equation}\label{x_l_def2}
X_l := \mathbb{E}({U}_n|{\cal G}_l) - \mathbb{E} ({U}_n|{\cal G}_{l-1}).
\end{equation}
We then have \(\sum_{l=1}^{N}X_l = U_n - \mathbb{E}U_n\) and
using the martingale property as in~(\ref{y_m_est}), we have
\[\mathbb{E}\left(\sum_{l=1}^{N}X_l\right)^2 = \sum_{l=1}^{N}\mathbb{E}X_l^2\] and so
\begin{equation}\label{prelim_est_t_n3}
\mathbb{E}(U_n-\mathbb{E}U_n)^2  = \sum_{l=1}^{N}\mathbb{E}X_l^2.
\end{equation}

As in (\ref{eq1}), we have the following estimate regarding \(X^2_l.\) We have
\begin{equation}
\mathbb{E}(X_l^2|{\cal G}_{l-1}) \leq C_2 \mathbb{P}(e_l \in \gamma_n | {\cal G}_{l-1}) \text{ a.s. }\label{eq2}
\end{equation}
for some positive constant \(C_2\) not depending on~\(l.\)
Using the above and proceeding as in the derivation of (\ref{fin_est_t_n}), we have
\[\mathbb{E}({U}_n -\mathbb{E}{U}_{n})^2 \leq C_2\mathbb{E}\left(\#\gamma_n\right), \] where \(\#\gamma_n\) denotes the number
of edges in the geodesic \(\pi_n.\) Substituting the estimate (\ref{len_pi_n2}) of Lemma~\ref{prop_pi_hat}, we have (\ref{eun_est2}).

To see (\ref{etn_est2}), we proceed as in the case of (\ref{etn_est}). As in (\ref{tn_un_est}), we have
\begin{equation}\label{tn_un_est2}
var(T_n) \leq 2 var(T_n - U_n)  + 2 var(U_n) \leq 2 var(T_n-U_n) + 2C_1n
\end{equation}
for some constant \(C_1 > 0\) and all \(n \geq 1,\) where the final estimate follows using (\ref{eun_est2}).
Set \(\epsilon = 3.\) For a fixed \(\eta > 0,\) there exists \(N_0 = N_0(\eta) \geq 1\) so that
\begin{equation}\label{un_est2}
var\left(T_n - U_n\right) \leq 2n^2\eta
\end{equation}
for all \(n \geq N_0(\eta).\)
Thus \[var(T_n) \leq 2C_1n  + n^2 \eta \leq 2n^2 \eta \] for all \(n \geq N_1(\eta).\) Since \(\eta > 0\) is arbitrary, this proves~(\ref{etn_est2}) of the Lemma.

%CHEK ALL DETAILS +eTC...

\emph{Proof of (\ref{un_est2})}: Analogous to the proof of (\ref{un_est}), we let
\begin{equation}\nonumber
H_n = \left\{\pi_n \subseteq B_{8\mu \beta_1^{-1} n^{1+\epsilon}}\right\}
\end{equation}
be the event that every edge of the geodesic~\(\pi_n\) for the passage time \(T_n\) is contained in the box~\(B_{8\mu \beta_1^{-1} n^{1+\epsilon}}.\)
From the definition of the passage times \(U_n\) and \(T_n\) we have that
\begin{equation}\label{tu_eq2}
U_n\ind(H_n) = T_n\ind(H_n).
\end{equation}

Letting \(Z_n = T_n - U_n,\) we have
\begin{equation}\label{sec_term_est2}
var({Z}_n) \leq \mathbb{E}Z_n^2 = \mathbb{E}Z_n^2\ind(H_n^c).
\end{equation}
We have that \[|Z_n| \leq T_n +U_n \leq 2\sum_{i=1}^{n}t(f_i) =: Q_n\] since the terms~\(U_n\) and~\(T_n\) are each no more than the passage time of the straight line with endvertices origin and~\((n,0,\ldots,0).\) We recall that for \(i \geq 1,\) the edge \(f_i\) is the edge joining \((i-1,0,\ldots,0)\) and \((i,0,\ldots,0).\) We therefore have from (\ref{sec_term_est2}) that
\begin{equation}\label{eq_zn1}
var(Z_n) \leq \mathbb{E}Z_n^2 \ind (H_n^c)  \leq \mathbb{E}Q_n^2 \ind (H_n^c) \leq 4n\sum_{i=1}^{n}\mathbb{E}t^2(f_i)\ind(H_n^c).
\end{equation}
For the final estimate, we use \((\sum_{i=1}^{k} a_i)^2 \leq k\sum_{i=1}^{k} a_i^2\) for positive \(a_i\) and have that \(Q_n^2 \leq 4n \sum_{i=1}^{n}t^2(f_i).\)
We evaluate each term in the above sum separately.

For a fixed \(1 \leq i \leq n,\) we have
\begin{equation}\label{eq_zn3}
\mathbb{E}t^2(f_i)\ind (H_n^c)  = \mathbb{E}t^2(f_i)\ind(H_n^c \cap \{t(f_i) \geq n^{\theta}\}) + \mathbb{E}t^2(f_i)\ind(H_n^c \cap \{t(f_i) < n^{\theta}\})
\end{equation}
for some constant \(\theta > 0\) to be determined later.
The second term is bounded above by
\begin{equation}\label{eq_znt1}
n^{2\theta}\mathbb{P}(H_n^c \cap \{t(f_i) < n^{\theta}\}) \leq n^{2\theta} \mathbb{P}(H_n^c) \leq n^{2\theta} \frac{C}{n^{1+2\epsilon}} = \frac{C}{n^{2\epsilon-2\theta+1}}
\end{equation}
where the final estimate is obtained using (\ref{box_pi}) of Lemma~\ref{box_est}.

Setting \(\theta = 1,\) the first term in (\ref{eq_zn3}) is bounded above by
\begin{equation}\label{eq_znt2}
\mathbb{E}t^2(f_i)\ind(t(f_i) \geq n) = \mathbb{E}t^2(f_i)\ind(t^2(f_i) \geq n^{2}) < \eta
\end{equation}
for all \(n \geq N_0.\) Here \(N_0 = N_0(\eta)\) depends only on \(\eta > 0\) and not on the choice of~\(i.\) The final estimate is true by the uniform integrability condition~\((ii)(a)\) of Section~\ref{intro}. Adding (\ref{eq_znt1}) and (\ref{eq_znt2}) and using (\ref{eq_zn3}) gives
\begin{equation}\label{h_nc_temp}
\mathbb{E}t^2(f_i)\ind (H_n^c)  \leq \frac{C}{n^{2\epsilon-1}} + \eta
\end{equation}
for all \(1 \leq i \leq n\) and all \(n \geq N_0.\)

Substituting the above estimate into (\ref{eq_zn1}) and setting \(\epsilon = 3,\) we have
\begin{equation}\nonumber
var(Z_n) \leq  4\frac{n^2C}{n^{2\epsilon-1}}  + 4n^2\eta = 4\frac{C}{n^{2\epsilon-3}} + 4n^2\eta \leq 4\frac{C}{n^3} + 4n^2\eta \leq 5n^2\eta.
\end{equation}
This proves (\ref{etn_est2}).\(\qed\)

To prove (\ref{etn_est3}), we proceed as in the case of (\ref{etn_est2}).
Set \(\epsilon = \frac{3}{p-2} + 1.\) We have that
\begin{equation}\label{un_est3}
var\left(T_n - U_n\right) \leq \frac{C}{n}
\end{equation}
for all \(n \geq 1\) and some constant \(C >0.\)
Substituting into (\ref{tn_un_est2}), we have that \[var(T_n) \leq 2C_1n  + \frac{C}{n} \leq 3C_1n\] for all \(n\) large and some constant \(C_1 >0.\) This proves~(\ref{etn_est3}) of the Lemma.
\emph{Proof of (\ref{un_est3})}: The proof proceeds as in (\ref{un_est2}) until (\ref{eq_znt2}). In particular (\ref{eq_znt1}) holds and instead of (\ref{eq_znt2}), we have
\begin{equation}\label{eq_znt3}
\mathbb{E}t^2(f_i)\ind(t(f_i) \geq n^{\theta}) = \mathbb{E}t^2(f_i)\ind(t(f_i) \geq n^{\theta}) \leq \left(\mathbb{E}t^p(f_i)\right)^\frac{2}{p}\left(\mathbb{P}(t(f_i) \geq n^{\theta})\right)^{1-\frac{2}{p}}
\end{equation}
where the final estimate follows from the Cauchy-Schwarz inequality.
Using the moment condition \((ii)(b),\) we have that \(\left(\mathbb{E}t^p(f_i)\right)^\frac{2}{p} \leq C_1\) for some constant \(C_1 > 0.\)
Also
\[\left(\mathbb{P}(t(f_i) \geq n^{\theta})\right)^{1-\frac{2}{p}} \leq \left(\frac{\mathbb{E}t^p(f_i)}{n^{p\theta}}\right)^{1-\frac{2}{p}} \leq \frac{C}{n^{\theta(p-2)}}\] for some constant \(C > 0.\) The first estimate above follows using Markov inequality and the final estimate follows from the moment condition \((ii)(b)\) of Section~\ref{intro}.

Setting \(\theta = \frac{3}{p-2}\) we have from (\ref{eq_znt3}) that \(\mathbb{E}t^2(f_i)\ind(t(f_i) \geq n^{\theta}) \leq \frac{C_2}{n^{3}}\) for some constant \(C_2 >0\) and setting \(\epsilon = \theta+1\) we have from (\ref{eq_znt1}) that
\[\mathbb{E}t^{2}(f_i)\ind(H_n^c \cap \{t(f_i) < n^{\theta}\}) \leq \frac{C}{n^{3}}\] for some constant \(C > 0.\) Thus
\[\mathbb{E}t^2(f_i)\ind (H_n^c)  \leq \frac{C}{n^{3}} + \frac{C_2}{n^{3}} = \frac{C_3}{n^3}\]
for all \(1 \leq i \leq n\) and all \(n \geq N_0\) large.
Substituting the above estimate into~(\ref{eq_zn1}) we have
\(var(Z_n) \leq  4\frac{n^2C_3}{n^{3}}  = 4\frac{C_3}{n}.\) This proves (\ref{un_est3}).\(\qed\)

\setcounter{equation}{0}
\renewcommand\theequation{\thesection.\arabic{equation}}
\section{Proofs of Theorem~\ref{thm1} and Corollary~\ref{cor1}}\label{pf_main_th}
\emph{Proof of (\ref{conv_tn_main}) of Theorem~\ref{thm1}}: We first prove a.s. convergence using a subsequence argument. From Proposition~\ref{lemma_v_n}, it suffices to prove that \(\frac{S_n}{n}\) converges to zero a.s. as \(n \rightarrow \infty,\) where  \[S_n = \hat{T}^{(n)}_n - \mathbb{E}\hat{T}^{(n)}_n.\] To prove that, we use a subsequence argument as follows. From Lemma~\ref{eun}, we have that \(\mathbb{E}S_n^2 \leq C n\) for some constant \(C > 0.\) Thus for a fixed \(\epsilon  > 0,\) we have that \[\sum_{n \geq 1} \mathbb{P}(|S_{n^2}| > n^2 \epsilon) \leq \sum_{n \geq 1}\frac{\mathbb{E}S^2_{n^2}}{\epsilon^2 n^4} \leq \sum_{n \geq 1} \frac{C}{\epsilon^2 n^{2}} < \infty.\] Since this is true for all \(\epsilon > 0,\) we have by Borel-Cantelli Lemma that
\begin{equation}\label{s_n_2_conv}
\frac{S_{n^2}}{n^2} \longrightarrow 0 \text{ a.s. }
\end{equation} as \(n \rightarrow \infty.\)

We now set
\begin{equation}\label{d_n_def}
D_{n^2} = \max_{n^2 \leq k < (n+1)^2} |S_k-S_{n^2}|
\end{equation}
and estimate \(D_{n^2}\) as follows. For \(n^2 \leq k < (n+1)^2,\) we write
\begin{eqnarray}
|S_k - S_{n^2}| &\leq& |\hat{T}^{(k)}_k- \hat{T}^{(n^2)}_{n^2}| + \mathbb{E}|\hat{T}^{(k)}_k - \hat{T}^{(n^2)}_{n^2}| \nonumber\\
&\leq& |\hat{T}^{(k)}_k- \hat{T}^{(k)}_{n^2}|  + |\hat{T}^{(k)}_{n^2}- \hat{T}^{(n^2)}_{n^2}| \nonumber\\
&&\;\;\;\;\;\;\;\;\;\;\; +\;\mathbb{E}|\hat{T}^{(k)}_k- \hat{T}^{(k)}_{n^2}|  + \mathbb{E}|\hat{T}^{(k)}_{n^2}- \hat{T}^{(n^2)}_{n^2}|.\;\;\;\;\label{eq_s2}
\end{eqnarray}
From (\ref{sub_add}) in Lemma~\ref{sec_mom} we have \[|\hat{T}^{(k)}_{k} - \hat{T}^{(k)}_{n^2}| \leq k^{\alpha}(k-n^2) \leq (n+1)^{2\alpha}((n+1)^2 - n^2) \leq C_1 n^{1+2\alpha}\]  for some constant \(C_1  >0.\) The second estimate holds since \(k < (n+1)^2.\) Substituting the above estimate into (\ref{eq_s2}), we obtain that
\begin{eqnarray}
|S_k - S_{n^2}| \leq 2C_1 n^{1+2\alpha}  + |\hat{T}^{(k)}_{n^2}- \hat{T}^{(n^2)}_{n^2}| + \mathbb{E}|\hat{T}^{(k)}_{n^2}- \hat{T}^{(n^2)}_{n^2}|.\;\;\;\;\label{eq_s3}
\end{eqnarray}

To estimate the remaining terms, we use estimate (\ref{tk1k2}) of Lemma~\ref{sec_mom} to obtain that \[0 \leq \hat{T}^{(k)}_{n^2} - \hat{T}^{(n^2)}_{n^2} \leq \hat{T}^{((n+1)^2)}_{n^2} - \hat{T}^{(n^2)}_{n^2}  =: I_{n^2}\] since \(n^2 \leq k < (n+1)^2.\) Thus letting \(D_{n^2}\) as in (\ref{d_n_def}), we have that  \[\frac{D_{n^2}}{n^2} \leq \frac{2C_1}{n^{1-2\alpha}} + \frac{I_{n^2}}{n^2} + \frac{\mathbb{E}I_{n^2}}{n^2}.\] We claim that \(\frac{I_{n^2}}{n^2} \rightarrow 0\) a.s. and that \(\frac{I_{n^2}}{n^2}\) is uniformly integrable. Assuming the claims for the moment, we then have \(\frac{\mathbb{E}I_{n^2}}{n^2} \rightarrow 0\) and since  \(\alpha < \frac{1}{2}\) (see (\ref{eps_def})), we get that  \(\frac{D_{n^2}}{n^2} \longrightarrow 0\) a.s. as \(n \rightarrow \infty.\) Also for \(n^2 \leq k < (n+1)^2, \) we have that \[\frac{|S_k|}{k}  \leq \frac{|S_k-S_{n^2}|}{k} + \frac{|S_{n^2}|}{k} \leq \frac{|S_k-S_{n^2}|}{n^2} + \frac{|S_{n^2}|}{n^2} \leq \frac{D_{n^2}}{n^2} + \frac{|S_{n^2}|}{n^2}.\] This proves that the original sequence \(\frac{S_k}{k} \rightarrow 0\) a.s. as \(k \rightarrow \infty.\)

To prove the two claims regarding \(I_{n^2},\) we use the fact that \[\hat{T}^{((n+1)^2)}_{n^2}\ind(W_n) = \hat{T}^{(n^2)}_{n^2}\ind(W_n)\] where \(W_n = \{T^{(n^2)}_{n^2} = T^{((n+1)^2)}_{n^2}\}\) is the event defined in (\ref{v_n_def}). From (\ref{w_n_est_tot}) of Lemma~\ref{lemma_v_n} and Borel-Cantelli lemma, we have that \(\mathbb{P}(\liminf_n W_n) = 1.\) Since \(I_{n^2} = I_{n^2}\ind(W_n^c),\) we get that \(\frac{I_{n^2}}{n^2} \rightarrow 0\) a.s. as \(n \rightarrow \infty.\)

To prove the uniform integrability of \(\frac{I_{n^2}}{n^2},\) we note that \[0 \leq I_{n^2} \leq \hat{T}^{((n+1)^2)}_{n^2} \leq \sum_{i=1}^{n^2}t^{((n+1)^2)}(f_i) \leq \sum_{i=1}^{n^2}t(f_i) \] where as before \(f_i\) denotes the edge from \((i-1,0,\ldots,0)\) to \((i,0,\ldots,0).\) From~(\ref{unif_int}) we therefore have \[\mathbb{E}I^2_{n^2} \leq \mathbb{E}\left(\sum_{i=1}^{n^2}t(f_i) \right)^2 \leq C n^4\] for some constant \(C > 0\) and so \(\mathbb{E}\left(\frac{I_{n^2}}{n^2}\right)^2 \leq C\) for all \(n \geq 1.\) This implies that \(\frac{I_{n^2}}{n^2}\) is uniformly integrable as desired. This proves the a.s. convergence in (\ref{conv_tn_main}).

To prove convergence in \(L^1,\) it is enough to see that \(\frac{T_n}{n}\) is uniformly integrable which follows from estimate (\ref{unif_int}) of Lemma~\ref{sec_mom}. Suppose now that the integrability condition \((ii)(a)\) holds, then we have from Lemma~\ref{eun2} that \(\frac{1}{n}(T_n-\mathbb{E}T_n)\) converges to zero in \(L^2.\)~\(\qed\)

%(//TO WRT CRFLLY...+ETC...\(\qed\)

%INCLUDE IN THEREOM 1...PRKVMM+ETC...

\emph{Proof of (\ref{tmin}) of Theorem~\ref{thm1}}: For the upper bound, we use (\ref{conv_tn_main}) and obtain a.s. that \(\limsup_n \frac{T_n}{n} = \limsup_n\frac{\mathbb{E}T_n}{n}.\) But \(T_n \leq \sum_{i=1}^{n} t(f_i)\) since the second term is the passage time of the straight line with endvertices origin and \((n,0,\ldots,0).\) Thus \(\mathbb{E}T_n \leq \sum_{i=1}^{n} \mathbb{E}t(f_i) \leq \mu n\) where \(\mu = \sup_i \mathbb{E}t(e_i) \in (0,\infty)\) is as in~(\ref{mu_def}). Thus \(\limsup_n \frac{\mathbb{E}T_n}{n} \leq \mu < \infty.\) This proves the upper bound in~(\ref{tmin}).

For the lower bound in (\ref{tmin}), we again use (\ref{conv_tn_main}) and obtain a.s. that \(\liminf \frac{T_n}{n} = \liminf \frac{\mathbb{E}T_n}{n}.\) It therefore suffices to see that the second term is positive. We proceed as follows. Suppose that the complement of the event~\(E_k\) defined in~(\ref{e_k_def}) occurs for \(k = k(n) = [\beta_1 (8\mu)^{-1}n].\) Here \(\beta_1 \in (0, \mu)\) is the constant in Lemma~\ref{t_pi_est} and~\([x]\) denotes the largest integer less than or equal to \(x.\) By our choice of \(k\) we have
\begin{equation}\label{eight_mu}
\frac{n}{2} \leq 8\mu \beta_1^{-1}k \leq n
\end{equation}
for all large \(n.\) We then have that every path containing \(r \geq 8\mu\beta_1^{-1} k\) edges has passage time at least \(\beta_1r \geq 8\mu k \geq \frac{\beta_1}{2}n,\) where the final estimate follows from (\ref{eight_mu}). Again using (\ref{eight_mu}), we have that every path containing at least~\(n\) edges has passage time at least \(\frac{\beta_1}{2}n.\)

Let \(Z_0\) be the null set in Proposition~\ref{geo_prop} so that for all \(\omega \in Z_0^c,\) finite paths attain the minimum passage time \(T_n.\) Fix \(\omega \in Z_0^c\) and let \(\pi_n\) be the path whose passage time \(T(\pi_n) = T_n\) (see (\ref{t_pi_def})). The path~\(\pi_n\) has at least~\(n\) edges and so if \(Z_0^c \cap E^c_k\) occurs, then by the discussion in the previous paragraph, we have \(T_n \geq \frac{\beta_1}{2}n.\) Thus
\begin{eqnarray}
\mathbb{E}T_n \geq \mathbb{E}T_n \ind(Z_0^c \cap E_k^c) \geq \frac{\beta_1}{2}n\mathbb{P}(Z_0^c \cap E^c_{k}) = \frac{\beta_1}{2}n\mathbb{P}(E^c_{k}) \geq \frac{\beta_1}{2}n(1-C_2e^{-\beta_2 k}), \nonumber
\end{eqnarray}
where \(C_2,\beta_2 > 0\) are as in (\ref{a_0k}). Using (\ref{eight_mu}) we have
\[1-C_2e^{-\beta_2 k} \geq 1-C_2e^{-\beta_3 n} \geq \frac{1}{2}\] for some constant \(\beta_3 > 0\) and for all \(n\) large. Thus \(\mathbb{E}T_n \geq \frac{\beta_1}{4}n\) for all \(n\) large and thus \(\liminf_n \frac{\mathbb{E}T_n}{n} \geq \frac{\beta_1}{4} > 0.\) This proves the lower bound in~(\ref{tmin}). \(\qed\)

\emph{Proof of Corollary~\ref{cor1}}: We show that \(\frac{\mathbb{E}T_n}{n} \rightarrow \mu_F\) for some constant \(\mu_F>0.\) First, using the version of subadditivity (\ref{sub_add}) for \(T_n\) we have
\begin{equation}\label{sub_add_exp}
\mathbb{E}T_{n+m} \leq \mathbb{E}T_n + \mathbb{E}T_{n,m+n} = \mathbb{E}T_n + \mathbb{E}T_m,
\end{equation}
where the final equality follows from translational invariance. As in the proof of (\ref{sub_add}), the term~\(T_{n,m+n}\) is the minimum passage time between the vertices \((n,0,\ldots,0)\) and \((m+n,0,\ldots,0).\) From (\ref{sub_add_exp}) and Fekete's Lemma we have that
\begin{equation}\label{lim_tn}
\lim_n\frac{\mathbb{E}T_n}{n} = \inf_{n \geq 1} \frac{\mathbb{E}T_n}{n} =: \mu_F.
\end{equation}
From (\ref{tmin}) of Theorem~\ref{conv_tn_main}, we have that \(\mu_F > 0.\) \(\qed\)

%\geq 4\mu n\] for all \(n\) large,  by (\ref{a_0k}). \(\qed\)

%\section{Future Work}\label{future}
%We briefly mention a couple of issues for possible future work. One possible question to consider is what further conditions need to be imposed on the passage times (in the non i.i.d. case) for convergence of the non-centred random variable \(\frac{T_n}{n}?\) As seen in Section~\ref{intro}, it is not necessarily always true that \(\frac{T_n}{n}\) converges.

%Secondly, what conditions would be needed to obtain in shape-type results in \(\mathbb{R}^2\) in this nonhomogenous scenario. More formally, we first extend the definition of \(T_x\) for all \(x \in \mathbb{R}^2\) by letting it be the minimum of the passage times to each of the four corners in \(\mathbb{Z}^2,\) of the \(1 \times 1\) square containing \(x.\) Letting \(A_t = \{x \in \mathbb{R}^2 : T_x \leq t\}\) we would like to know if there exists a function \(\varphi : \mathbb{R}^2 \rightarrow \mathbb{R}\) such that \[\mathbb{P}\left(\left\{ L_{\epsilon} \subseteq \frac{A_t}{t} \subseteq U_{\epsilon} \text{ for all large t}\right\}\right) = 1,\] where \(L_{\epsilon} = \{x : \varphi(x) \leq 1 - \epsilon\}\) and \(U_{\epsilon} = \{x : \varphi(x) \leq 1+\epsilon\}.\) For further studies on shape results in homogenous first passage percolation we refer to Cox and Durrett (1981), Kesten (1993) and references therein. %and if so, can we employ the scaling methods of Cox et. al. (1993)?

\section*{Acknowledgements}
I thank Professors Rahul Roy, Thomas Mountford, Siva Athreya and Federico Camia for crucial comments and for fellowships. I also thank NBHM and ISF-UGC for my fellowships.

%\subsection*{Left-right crossings in \(\mathbb{Z}^2\)}
%Let \(R\) be any \(\frac{m_1 r_n}{4}  \times \frac{m_2 r_n}{4}\) rectangle that is contained in \(S\) that contains exactly \(m_1m_2\) of the squares \(S_i.\) There exists a corresponding \(m_1 \times m_2\) rectangle, \(\tilde{R} \subset \mathbb{Z}^2\) (see Figure~\ref{rect_comp2} for an illustration). We define an open left-right crossing~\cite{Grimmett} in \(\tilde{R}\) to be any open connected path of sites from the leftmost set of vertices to the rightmost in \(\tilde{R}.\) In Figure~\ref{rect_comp2}, the open sites are circled in the rectangle \(\tilde{R}.\) The thin dotted line forms an open left-right crossing of \(\tilde{R}\) from \(a\) to \(b.\) For all sufficiently large \(n,\) we claim that the event that a left-right crossing occurs in \(\tilde{R}\) has probability at least
%\begin{equation}\label{lr_prob}
%1 - 2m_2(4(1-p))^{m_1}
%1 - \frac{m_2}{n^{m_1\delta_1}}
%\end{equation}
%for some \(\delta_1 > 0.\) %Assuming that (\ref{lr_prob}) holds for the moment, we complete the rest of the proof.

%erere

\setcounter{equation}{0} \setcounter{Lemma}{0} \renewcommand{\theLemma}{II.%
\arabic{Lemma}} \renewcommand{\theequation}{II.\arabic{equation}} %
\setlength{\parindent}{0pt}

%\section*{Appendix~II: Proof of (\ref{p_t_y_ineq1}), (\ref{e_2_m_c1}) and (\ref{e_2_m_c})}\label{app1}
%\begin{figure}
%\centering
%\includegraphics[width=3.5in]{cap_scheme.eps}
%\caption{Modified coding scheme of~\cite{Goldsmith}.} \label{fig_eta_res}
%\end{figure}

%\begin{figure}
%\centering
%\includegraphics[width=4in]{cap_vs_se.eps}
%\caption{Capacity as a function of network population \(n\) for various values of \(\sigma_e.\)} \label{cap_se}
%\end{figure}

%\bibliographystyle{IEEEtran}
%\bibliography{IEEEabrv,ref_mob}

%\begin{biography}[{\includegraphics[width=1in,height=1.25in,clip,keepaspectratio]{mshell}}]{Michael Shell}
% where an .eps filename suffix will be assumed under latex, and a .pdf suffix
% will be assumed for pdflatex; or if you just want to reserve a space for
% a photo:


%\vfill



\bibliographystyle{plain}

\end{document}